\documentclass[12pt, draft]{article}
\usepackage{amssymb}
\usepackage{amsmath}
\usepackage{array}
\usepackage[T2A]{fontenc}
\usepackage[cp1251]{inputenc}
\usepackage[russian]{babel}
\usepackage[pdftex]{graphics}

\renewcommand{\S}{\mathhexbox278}

\DeclareMathOperator{\RRe}{Re} \DeclareMathOperator{\IIm}{Im}

\DeclareMathOperator{\vth}{\vartheta}
\DeclareMathOperator{\vep}{\varepsilon}
\DeclareMathOperator{\vk}{\varkappa}

\renewcommand{\le}{\operatorname{\leqslant}}
\renewcommand{\ge}{\operatorname{\geqslant}}

\multlinegap=0em \DeclareFontFamily{T1}{msb}{}
\DeclareFontShape{T1}{msb}{m}{ol}{<5> <6> <7> <8> <9> gen * msbm
<10> <10.95> <12> <14.4> <17.28> <20.74> <24.88> msbm10}{}
\DeclareSymbolFont{AMSb}{T1}{msb}{m}{ol} \multlinegap=0em
\begin{document}

\begin{center}
{\rmfamily\bfseries\normalsize Gram's Law and the Argument of the
Riemann Zeta Function\footnote{This research was supported by the
Russian Fund of Fundamental Research (grant no. 11-01-00759-a).}}
\end{center}

\begin{center}
M.A. Korolev

Steklov Mathematical Institute

Gubkina str., 8,

119991, Moscow, Russia

\texttt{hardy\_ramanujan@mail.ru}, \texttt{korolevma@mi.ras.ru}
\end{center}

\vspace{0.5cm}

\fontsize{11}{12pt}\selectfont

\textbf{Annotation.} Some new statements concerning the behavior of
the argument of Riemann zeta function at the Gram points are proved.
We apply these statements to the proof of Selberg's formulas
connected with Gram's Law.


\fontsize{12}{15pt}\selectfont

\begin{center}
\textbf{\S 1. Introduction}
\end{center}

Since the notion `Gram's Law' has different sense in different
papers, it seems re\-a\-son\-ab\-le to begin this paper with short
survey. This survey contains basic results concerning the peculiar
phenomenon observed by J\"{o}rgen Pedersen Gram \cite{Gram_1903} in
1903.

For the following, we need several definitions. Suppose $t>0$ and
let $\vth(t)$ be an increment of any fixed continuous branch of the
argument of the function $\pi^{-s/2}\Gamma\bigl(\tfrac{s}{2}\bigr)$
as $t$ varies along the segment connecting the points $s =
\tfrac{1}{2}$ and $s = \tfrac{1}{2}+it$. Then Hardy's function
$Z(\tau) = e^{i\vth(\tau)}\zeta\bigl(\tfrac{1}{2}+i\tau\bigr)$ is
real for real $\tau$ and it's real zeros coincide with the ordinates
of zeros of $\zeta(s)$ lying on the critical line $\RRe s =
\tfrac{1}{2}$. Further, if $t$ is not an ordinate of a zero of
$\zeta(s)$ then the function $S(t) =
\pi^{-1}\arg\zeta\bigl(\tfrac{1}{2}+it\bigr)$ is defined as an
increment of any continuous branch of $\pi^{-1}\arg\zeta(s)$ along
the polygonal arc connecting the points $2, 2+it$ and
$\tfrac{1}{2}+it$. In the opposite case, $S(t)$ is defined by the
relation
\[
S(t)\,=\,\tfrac{1}{2}\lim_{h\to 0}\bigl(S(t+h)+S(t-h)\bigr).
\]
Let $N(t)$ be a number of zeros of $\zeta(s)$ in the rectangle
$0<\IIm s \le t$, $0\le \RRe s \le 1$ counted with multiplicities.
At the points of discontinuity $N(t)$ is defined as follows:
\[
N(t)\,=\,\tfrac{1}{2}\lim_{h\to 0}\bigl(N(t+h)+N(t-h)\bigr).
\]
The equation
\begin{equation}
N(t)\,=\,\frac{1}{\pi}\,\vth(t)\,+\,1\,+\,S(t)\label{Lab01}
\end{equation}
holds true for any $t>0$ and is called as Riemann - von Mangoldt
formula.

By $\varrho_{n} = \beta_{n}+i\gamma_{n}$ we denote complex zeros of
$\zeta(s)$ lying in the upper half-plane and ordered as follows:
$0<\gamma_{1}<\gamma_{2}<\ldots \le \gamma_{n}\le \gamma_{n+1}\le
\ldots$. Finally, let $c_{n}$ be a positive zeros of $Z(t)$ numbered
in ascending order and counted with their multiplicities.

Though first three positive ordinates of zeros of $\zeta(s)$ had
been counted already by Riemann, this fact was declared by
C.L.Siegel \cite{Siegel_1932} in 1932 only. It seems that first
mathematical publication devoted to the calculation of zeta zeros
belongs to Gram \cite{Gram_1895} (1895). He established that
$\gamma_{1}\,=\,14.135, \gamma_{2}\,=\,20.82, \gamma_{3}\,=\,25.1$,
but his method was too laborious and unfit for the detecting of
higher zeros. In 1902, Gram invented a more acceptable method of
detecting of the zeros of $\zeta(s)$.

The key idea of this method is the following. Let $A(t)$ and $B(t)$
be a real and imaginary parts of $\zeta\bigl(\tfrac{1}{2}+it\bigr)$
correspondingly. Then
\[
\zeta\bigl(\tfrac{1}{2}+it\bigr)\,=\,e^{-i\vth(t)}Z(t)\,=\,Z(t)(\cos{\vth(t)}\,-\,i\sin{\vth(t)}),
\]
and hence $A(t) = Z(t)\cos{\vth(t)}$, $B(t) = -Z(t)\sin{\vth(t)}$.
Now let us consider the real zeros of $B(t)$. These zeros are of two
types. The zeros of the first type are the ordinates $\gamma_{n}$ of
zeros of $\zeta(s)$ lying on the critical line, and the zeros of the
second type are the roots of the equation $\sin{\vth(t)} = 0$. Using
the Stirling's formula in the form
\[
\vth(t)\,=\,\frac{t}{2}\ln\frac{t}{2\pi}\,-\,\frac{t}{2}\,-\,\frac{\pi}{8}\,+\,O\Bigl(\frac{1}{t}\Bigr),
\]
and considering the values $t>7$, it's possible to show that the
roots of the above equation generate the unbounded monotonic
sequence: $t_{0} = 9.6669\ldots$, $t_{1} = 17.8456\ldots$, $t_{2} =
23.1703\ldots$, $t_{3} = 27.6702\ldots$\,, $\ldots$\,. Here $t_{n}$
denotes $n$th Gram point i.e. the unique solution of the equation
$\vth(t_{n}) = (n-1)\pi$. Thus, the value
$\zeta\bigl(\tfrac{1}{2}+it_{n}\bigr)$ is real and
\[
\zeta\bigl(\tfrac{1}{2}+it_{n}\bigr)\,=\,A(t_{n})\,=\,Z(t_{n})\cos{\pi(n-1)}\,=\,(-1)^{n-1}Z(t_{n}).
\]

Suppose now that $A(t)$ has the same sign at the points $t_{n-1}$
and $t_{n}$ for some $n$. Then the values $Z(t_{n-1})$ and
$Z(t_{n})$ are of opposite sign. Therefore, $Z(t)$ vanishes at the
odd number of points between $t_{n-1}$ and $t_{n}$.

Using Euler - Maclaurin summation formula, Gram established that
$A(t_{n})>0$ for $n = 1,2,\ldots, 15$, and proved that all zeros of
$\zeta(s)$ in the strip $0<t<66$ lie on the critical line. This
method allowed him also to find approximately the ordinates
$\gamma_{1}, \gamma_{2}, \ldots, \gamma_{15}$. Thus Gram established
that there is exactly one zero $c_{n}$ in each interval $G_{n} =
(t_{n-1}, t_{n}]$, $n = 1,2,\ldots, 15$ and, moreover,
$t_{n-1}<c_{n}<t_{n}$. However, he assumed also that this law is not
universal: `... \emph{the values} $A(t_{n})$ \emph{are positive for
all} $t_{n}$ \emph{lying between} $10$ \emph{and} $65$. \emph{It
seems that the function} $A(t)$ \emph{is positive for a most part
of} $t$ \emph{under consideration. Obviously}, \emph{the reason is
that the first term of the sum}
$\sum\limits_{1}^{n}n^{-1/2}\cos{(t\log{n})}$ \emph{leads to the
dominance of positive summands}. \emph{If it is so}, \emph{the
regularity in the relative location of} $c_{n}$ \emph{and} $t_{n}$
\emph{will hold true for some time for the roots} $c$ \emph{lying
closely to} $c_{15}$ \emph{until the equilibrium will set in}  (see
\cite{Gram_1903}; for simplicity, we use here the notations of the
present paper).

The words `Gram's law' appeared for the first time in Hutchinson's
paper \cite{Hutchinson_1925}. He used this notion to underline the
property of $c_{n}$ and $c_{n+1}$ to be separated by Gram's point
$t_{n}$. Hutchinson undertook more wide calculations of zeros of
$\zeta(s)$ in order to check the validity of Gram's assumption. He
found two first values of $n$ that do not satisfy Gram's law: $n =
127$ and $n = 136$. Namely, he established that
\[
t_{127}\,<\,\gamma_{127}\,<\,\gamma_{128}\,<\,t_{128},\quad\quad
t_{134}\,<\,\gamma_{135}\,<\,\gamma_{136}\,<\,t_{135}.
\]
Ten years later, E.C.Titchmarsh and L.J.Comrie
\cite{Titchmarsh_1935},\cite{Titchmarsh_1936} continued Hutchinson's
cal\-cu\-la\-ti\-ons using Brunsviga, National and Hollerith
machines. They found a lot of new ex\-cep\-ti\-ons from Gram's law,
but the proportion of these exceptions does not exceed
$4.5\%$\footnote{Titchmarsh mentioned 43 exceptions that he had
found during the calculation of first 1041 zeros of $Z(t)$ lying in
the interval $0<t\le 1468$. However, there are 1042 zeros of Hardy's
function and 1041 Gram's points $t_{n}$ between $t = 0$ and $t =
1468$, and there are 45 values of $n$ such that
$(-1)^{n-1}Z(t_{n})<0$.}. The paper \cite{Titchmarsh_1935} contains
also the first theoretical results concerning Gram's law. Thus,
Tit\-c\-h\-marsh proved that the inequality $A(t_{n}) =
(-1)^{n-1}Z(t_{n})>0$ fails for infinitely many $n$. Mo\-re\-over,
he proved that the sequence of fractions
$\tau_{n}\,=\,\frac{c_{n}-t_{n}}{t_{n+1}-t_{n}}$ is unbounded. The
last assertion means that there are infinitely many zeros $c_{n}$
lying outside the interval $G_{n}$.

Though the rule with infinitely many exceptions is not a law in a
rigorous sense,  the notion `Gram's law' is widely used now, but
sometimes in different senses. We will also use this notion for any
assertion concerning the relative location of ordinates of zeros of
$\zeta(s)$ and Gram points. In what follows, we give a sort of
`classification' of `Gram's laws'.

\vspace{0.2cm}

\textsc{Definition 1}. Gram's interval $G_{n} = (t_{n-1},t_{n}]$
satisfies to the Strict Gram's Law (SGL) iff $G_{n}$ contains a zero
$c_{n}$ of $Z(t)$.

\vspace{0.2cm}

This definition is close to that Hutchinson and Titchmarsh used. But
here we allow the coincidence of $t_{n}$ and a zero of Hardy's
function. The reason is that now very little is known about the
number of such coincidences (or noncoincidences). It seems that
$Z(t_{n})$ does not vanish for every $n$, i.e. $c_{m}\ne t_{n}$ for
any $m$ and $n$. We only know that $Z(t_{n})\ne 0$ for at least
$(4-o(1))N(\ln N)^{-1}$ values of $n$, $1\le n\le N$ (see
\cite{Kalpokas_Steuding_2011})\footnote{Selberg's theorem formulated
without proof in \cite{Selberg_1946b} and cited below implies that
$Z(t_{n})\ne 0$ for positive proportion of $n$. It's interesting to
note that the values $Z(t_{n})$ are very close to zero for some $n$.
Thus, the minima of $|Z(t_{n})|$ for $n\le 10^{5}$ and $n\le 10^{6}$
are equal to $1.238\cdot 10^{-5}$ ($n = 97\,281$) and to $8.908\cdot
10^{-8}$ ($n = 368\,383$) correspondingly.}. The unboundedness of
the fractions $\tau_{n}$ implies that SGL fails for infinitely many
cases. Unfortunately, it is still unknown whether the number of
cases when SGL holds true is finite or infinite.

The definition 1 contains a very heavy restriction. Namely, the
number of the interval $G_{n}$ and the number of zero $c$ belonging
to $G_{n}$ must be equal. The renunciation of this restriction leads
to the second version of Gram's law.

\vspace{0.2cm}

\textsc{Definition 2}. Gram's interval $G_{n}$ satisfies to Gram's
law (GL) iff $G_{n}$ contains exactly one (simple) zero of $Z(t)$.

\vspace{0.2cm}

It's possible to show that SGL and GL are not equivalent to each
other. The failure (validity) of one statement does not imply the
failure (validity) of another. For example, $G_{1},\ldots , G_{126}$
satisfy both SGL and GL; $G_{127}$ does not satisfy neither SGL, nor
GL; further, $G_{128}$ satisfies SGL, but does not satisfy GL;
finally, $G_{3359}$, $G_{3778}$, $G_{4542}$ satisfy GL, but do not
satisfy SGL.

The counting of zeros of $Z(t)$ in given interval $(a,b)$ usually
reduces to the determining of number of sign-changes of $Z(t)$ in
$(a,b)$. Therefore, this method allows to determine the evenness of
the number of zeros only. For example, the inequality $Z(a)Z(b)<0$
guarantees the existence of an odd number of zeros in $(a,b)$
counted with multiplicity. Therefore, it seems natural to consider
one more type of Gram's law.

\vspace{0.2cm}

\textsc{Definition 3}. Gram's interval $G_{n}$ satisfies to Weak
Gram's Law (WGL) iff $G_{n}$ contains an odd number of zeros of
$Z(t)$.

\vspace{0.2cm}

Obviously, GL implies WGL, but the opposite statement is incorrect.
Thus, if $n = 2147$ then $G_{n}$ contains exactly three zeros of
Hardy's function, namely $c_{n-1}, c_{n}$  and $c_{n+1}$. The
inequality $Z(t_{n-1})Z(t_{n})<0$ is a sufficient (but not
necessary) condition for WGL. Therefore, Titchmarsh's formula (see
\cite{Titchmarsh_1934})
\[
\sum\limits_{n\le N}Z(t_{n-1})Z(t_{n})\,\sim -\,2(\gamma + 1)N,
\]
implies that WGL holds true in infinitely many cases ($\gamma$
denotes Euler's constant).

The statement `SGL holds true for all $n\ge n_{0}$' implies the
boundedness of the fractions $\tau_{n}$ as $n\to +\infty$. The last
fact contradicts to some properties of $S(t)$ established by H.Bohr
and E.Landau \cite{Bohr_Landau_1913} in 1913.

In the middle of 40's, Selberg invented a new powerful method of
researching of $S(t)$ (see \cite{Selberg_1946a},
\cite{Selberg_1944}) and obtained a lot of very deep results
concerning the distribution of zeros of $\zeta(s)$. In particular,
in \cite{Selberg_1946b} he formulated (without proof) the following
theorem: there exist absolute constants $K$ and $N_{0}$ such that
for $N>N_{0}$, $1\le n\le N$, the numbers $Z(t_{n-1})$ and
$Z(t_{n})$ are of different sign in more than $KN$ cases, and of the
same sign in more than $KN$ cases. This theorem implies that both
WGL and GL fail for positive proportion of cases, and that WGL holds
true for positive proportion of cases.

Denote by $\nu_{k} = \nu_{k}(N)$ the number of Gram's intervals
$G_{n}$, $1\le n\le N$, that contain exactly $k$ ordinates of zeros
of $\zeta(s)$ (here we consider all the zeros in the critical strip,
but not only the zeros lying on the critical line). It's not
difficult to prove that Selberg's theorem implies the following
relations:
\begin{align}
& \nu_{0}\,+\,\nu_{2}\,+\,\nu_{4}\,+\,\ldots\,\ge\,KN, \label{Lab02}\\
& \nu_{1}\,+\,\nu_{3}\,+\,\nu_{5}\,+\,\ldots\,\ge\, KN,
\label{Lab03}
\end{align}

These inequalities are weaker than the original Selberg's assertion.
The reason is that Selberg's theorem deals with the ordinates in
open intervals $(t_{n-1},t_{n})$ and with non\,-vanishing of $Z(t)$
at the end\,-points of such intervals for positive proportion of
$n$.

Further, (\ref{Lab02}) implies weaker estimate
\begin{equation}
\nu_{0}\,+\,\nu_{2}\,+\,\nu_{3}\,+\,\nu_{4}\,+\ldots\,\ge
\,KN.\label{Lab04}
\end{equation}
This inequality shows that the positive proportion of Gram's
intervals contain an `ab\-nor\-mal' number (i.e. $\ne 1$) of
ordinates. Both (\ref{Lab02}) and (\ref{Lab04}) imply that WGL and
GL fail for positive proportion of $n$.

As far as author knows, the proof of Selberg's theorem or the proof
of (\ref{Lab02}), (\ref{Lab03}) were never published. The estimate
(\ref{Lab04}) was proved by T.S.~Trudgian \cite{Trudgian_2009} in
2009\footnote{It is necessary to note that the inequalities
$\nu_{0}\gg N$, $\nu_{2}+\nu_{3}+\nu_{4}+\ldots \gg N$ were
formulated without proof by A.~Fujii in \cite{Fujii_1987}}. He also
pointed out in \cite{Trudgian_2009} that (\ref{Lab04}) implies the
inequality
\begin{equation}
\nu_{0}\,>\,K_{1}N \label{Lab05}
\end{equation}
for any fixed $K_{1}$, $0<K_{1}<\tfrac{1}{2}K$, and for $N\ge
N_{0}(K_{1})$. Indeed, the following identities hold true:
\begin{align}
& 0\cdot\nu_{0}\,+\,1\cdot\nu_{1}\,+\,2\cdot\nu_{2}\,+\,3\cdot\nu_{3}\,+\ldots\,=\,N+S(t_{N}+0), \label{Lab06}\\
& \nu_{0}\,+\,\nu_{1}\,+\,\nu_{2}\,+\,\nu_{3}\,+\,\ldots\,=\, N.
\label{Lab07}
\end{align}
It's easy to see that (\ref{Lab06}) expresses the fact that the
number of zeros whose ordinates are positive and do not exceed
$t_{N}$, is equal to $N(t_{N}+0) = N+S(t_{N}+0)$, and (\ref{Lab07})
expresses the fact that the number of $G_{n}$ containing in
$(0,t_{N}]$ is equal to $N$. Subtracting (\ref{Lab07}) from
(\ref{Lab06}) we find:
\[
\nu_{0}\,=\,\nu_{2}\,+\,2\nu_{3}\,+\,3\nu_{4}\,+\,\ldots\,-\,S(t_{N}+0).
\]
Adding $\nu_{0}$ to both parts and using the classical estimate
$S(t) = O(\ln t)$ (see \cite{Mangoldt_1905}) we get:
\begin{multline*}
2\nu_{0}\,=\,\nu_{0}\,+\,\nu_{2}\,+\,2\nu_{3}\,+\,3\nu_{4}\,+\ldots\,\ge\\
\ge\,\nu_{0}\,+\,\nu_{2}\,+\,nu_{3}\,+\,\nu_{4}\,+\ldots\,+O(\ln
N)\,\ge KN\,+\,O(\ln N).
\end{multline*}

This proves (\ref{Lab05}). Similarly to (\ref{Lab02}) and
(\ref{Lab04}), the inequality (\ref{Lab05}) implies that both WGL
and GL fail for positive proportion of cases.

It's interesting to note the following. It is supposed that
$\nu_{1}(N)\ge cN$ or even $\nu_{1}(N)\sim c_{1}N$ as $N$ growths.
However, weaker relation $\nu_{1}(N)\to +\infty$ as $N\to +\infty$
is still unproved. Thus, we don't know whether the number of cases
when GL holds true is finite or not.

There are some reasons to think that Selberg interpreted Gram's Law
in a way different from Titchmarsh's one and different from SGL, GL
and WGL. In dealing with Gram's Law, Titchmarsh considered only the
real zeros of Hardy's function. We think that Selberg considered all
the zeros of $\zeta(s)$ in the critical strip. The basic arguments
that sustain this point of view will be introduced later. Now we
give here some preliminary notes.

Let $\gamma_{n}$ be an ordinate of any zero of $\zeta(s)$ in the
critical strip. Then we determine a unique integer $m = m(n)$ such
that $t_{m-1}\,<\,\gamma_{n}\,\le\,t_{m}$, and set $\Delta_{n} =
m-n$.

\vspace{0.2cm}

\textsc{Definition 4.} We say that Gram - Selberg's Phenomenon (GSP)
is observed for the ordinate $\gamma_{n}$ iff $\Delta_{n} = 0$, i.e.
iff $t_{n-1}<\gamma_{n}\le t_{n}$.

\vspace{0.2cm}

It seems very probable that the property of $\gamma_{n}$ to satisfy
the condition $\Delta_{n} = 0$ was called by Selberg as `Gram's
Law'.

The above result of Selberg implies that there is a positive
proportion of cases when GSP does not observed. However, it's
possible to say much more about GSP. Thus, Selberg established the
formulas
\begin{align}
& \sum\limits_{n\le
N}\Delta_{n}^{2k}\,=\,\frac{(2k)!}{k!}\,\frac{N}{(2\pi)^{2k}}\,(\ln\ln
N)^{k}\,+\,O\bigl(N(\ln\ln N)^{k-1/2}\bigr), \label{Lab08}\\
& \sum\limits_{n\le N}\Delta_{n}^{2k-1}\,=\,O\bigl(N(\ln\ln
N)^{k-1}\bigr), \label{Lab09}
\end{align}
where $k\ge 1$ is a fixed integer, and assumed that the inequalities
\[
\frac{1}{\Phi(n)}\,\sqrt{\ln\ln
n}\,<\,|\Delta_{n}|\,\le\,\Phi(n)\sqrt{\ln\ln n}
\]
hold true for `almost all' $n$. Here $\Phi(x)$ denotes any fixed
positive unbounded function. In particular, this assumption implies
that GSP does not observed in `almost all' cases.

It follows from the remark on p.355 of \cite{Selberg_1989} that
Selberg found a proof of his own assumption a long before 1989, but
he didn't publish it. For the reconstruction of such proof, see
author's papers \cite{Korolev_2010a},\cite{Korolev_2010b}.

Thus, we can't expect the occurrence of GSP in positive proportion
of cases. The reason is that the definition 4 impose a heavy
restriction on $\gamma_{n}$ (the ordinate should belong to the
Gram's interval with the same number). Thus, GSP is a very rare
phenomenon. It's natural to ask whether GSP occurs in infinitely
many cases or no. Since no such results were published, some
quantitative statements about the frequency of occurrence of GSP
seem to have some interest. Author intends to introduce them in a
future paper.

Now we place the results concerning Gram's law in the following
table.

\begin{center}
\begin{tabular}{|c|m{6.5cm}|m{6.5cm}|}
\hline & holds true & fails \\
\hline

SGL & it's unknown, whether the number of such cases is finite or no
& for infinitely many cases:

- Titchmarsh \cite{Titchmarsh_1935}, 1935; \\ \cline{1-1}
\cline{3-3}

GL & & for infinitely many cases:

- Titchmarsh \cite{Titchmarsh_1935}, 1935;

for positive proportion of cases:

- Selberg \cite{Selberg_1946b}, 1946;

- Fujii \cite{Fujii_1987}, 1987;

- Trudgian \cite{Trudgian_2009}, 2009. \\ \hline

WGL & for infinitely many cases:

- Titchmarsh \cite{Titchmarsh_1934}, 1934;

for positive proportion of cases:

- Selberg \cite{Selberg_1946b}, 1946;

& for positive proportion of cases:

- Selberg \cite{Selberg_1946b}, 1946;

- Fujii \cite{Fujii_1987}, 1987;

- Trudgian \cite{Trudgian_2009}, 2009.\\ \hline

GSP & for infinitely many cases & for infinitely many cases:

- Titchmarsh \cite{Titchmarsh_1935}, 1935;

for `almost all' cases:

- Selberg \cite{Selberg_1989}, 1989.\\
\hline
\end{tabular}
\end{center}

The present paper contains some new statements concerning the
behavior of the function $S(t)$ at Gram points. These statements are
applied to the proof of Selberg's formulas (\ref{Lab08}),
(\ref{Lab09}) and to other problems connected with Gram's law. The
paper is organized as follows.

First, \S 2 contains auxiliary assertions. In \S 3, the sum
\begin{equation}
\sum\limits_{N<n\le N+M}\bigl(S(t_{n+m}+0)\,-\,S(t_{n}+0)\bigr)^{2k}
\label{Lab10}
\end{equation}
is calculated. Here $k\ge 0, m>0$ are sufficiently large integers,
that may growth slowly with $N$ (Theorem 1). The correct bound (in
the sense of order of growth) for the sum (\ref{Lab10}) with $m = 1$
is also given here (Theorem 2). It is necessary to note that the
bounds of such type are contained in \cite{Fujii_1987}. But they
hold true only for a `long' interval of summation: $1\le n\le N$ or
$N<n\le N+M$, $M\asymp N$. The statements of present paper are valid
for a `short' interval, namely for the case $M\asymp
N^{\alpha+\vep}$, $\alpha = \tfrac{27}{82} =
\tfrac{1}{3}-\tfrac{1}{246}$.

Further, Theorem 3 in \S 3 gives a true order of magnitude of the
sum
\[
\sum\limits_{N<n\le N+M}\bigl|S(t_{n}+0)\,-\,S(t_{n-1}+0)\bigr|.
\]
This statement is based on Theorems 1 and 2 and plays the key role
in the proof of the inequalities
\[
\nu_{0}(N+M)\,-\,\nu_{0}(N)\,\gg\,M,\quad \sum\limits_{k\ge
2}\bigl(\nu_{k}(N+M)\,-\,\nu_{k}(N)\bigr)\,\gg\,M
\]
(Theorem 4). The analogues of these bounds for the case of the
`long' intervals of summation were formulated (without proof) for
the first time by Selberg in \cite{Selberg_1946b} as a corollary of
his theorem cited above.

In \S 4, a non-trivial bound for the alternating sum
\[
T_{k}\,=\,\sum\limits_{N<n\le
N+M}S^{k}(t_{n}+0)\bigl(S(t_{n}+0)\,-\,S(t_{n-1}+0)\bigr)
\]
is given (Theorem 5). This estimation leads to a new proof of
Selberg's formulas (\ref{Lab08}), (\ref{Lab09}) (see Theorem 6) and
of the assumption that $\Delta_{n}\ne 0$ for `almost all' $n$. The
proofs of these facts differ from those given in previous author's
paper \cite{Korolev_2010b}. They don't use the information about the
number of solutions of the inequalities $a<\Delta_{n}\le b$ with the
condition $N<n\le N+M$. It's possible that the below proof of
Theorem 6 is close to the original Selberg's proof.

Finally, in \S 6 we argue that Selberg considered all the complex
zeros of $\zeta(s)$ in dealing with the Gram's law in
\cite{Selberg_1946b}. This question leads us to a new equivalent of
`almost Riemann hypothesis' (see Theorem 7; `almost Riemann
hypothesis' claims that almost all complex zeros of $\zeta(s)$ lie
on the critical line).

Throughout the paper, $\vep$ denotes an arbitrary small positive
number, $0<\vep<10^{-3}$; $N_{0}(T)$ denotes the number of zeros of
$\zeta\bigl(\tfrac{1}{2}+it\bigr)$ for $0<t\le T$; $N\ge
N_{1}(\vep)>0$ is a sufficiently large integer; $L = \ln\ln N$, $M$
is an arbitrary integer with the conditions $N^{\alpha+\vep_{1}}\le
M\le N^{\alpha + \vep}$, $\alpha = \tfrac{27}{82}$, $\vep_{1} =
0.9\vep$; $\theta, \theta_{1}, \theta_{2}, \ldots$ are complex
numbers whose absolute values do not exceed 1 and which are,
generally speaking, different in different relations. In some cases
we use for brevity the notation $\Delta(n)$ for the value
$S(t_{n}+0)$.

\begin{center}
\textbf{\S 2. Auxiliary lemmas}
\end{center}

\textsc{Lemma 1}. \emph{The following relations hold true for any}
$x\ge 2$:
\[
\sum\limits_{p\le x}\frac{\ln p}{p}\,<\,\ln x, \quad\quad
\sum\limits_{p\le x}\frac{1}{p}\,=\,\ln\ln
x\,+\,c\,+\,\frac{\theta}{\ln^{2}x};
\]
\emph{here} $c = 0.26\ldots$ \emph{and} $-\tfrac{1}{2}<\theta<1$.

\vspace{0.2cm}

For a proof, see \cite{Rosser_Schoenfeld_1962}.

\vspace{0.2cm}

\textsc{Lemma 2}. \emph{Suppose that} $0<\kappa < \tfrac{1}{2}$,
$0<c<\tfrac{1}{2}-\kappa$, $\mu,\nu$ \emph{are integers such that}
$\mu, \nu \ge 0$, $\mu+\nu = 2k$, $k\ge 1$, $N\ge
\exp{(9\kappa^{-1})}$, $M\ge \exp{(3kc^{-1})}$, $y = M^{c/k}$.
\emph{Furthermore}, \emph{let} $p_{1},\ldots, p_{\nu}, q_{1},
\ldots, q_{\mu}$ \emph{take values of prime numbers in the interval}
$(1,y]$ \emph{and satisfy the condition} $p_{1}\ldots p_{\nu}\ne
q_{1}\ldots q_{\mu}$. \emph{Finally}, \emph{suppose that} $|a(p)|\le
\delta$ \emph{for} $p\le y$. \emph{Then the sum} $S$,
\[
S\,=\,\sum\limits_{N<n\le N+M}\sum\limits_{\substack{p_{1},\ldots,
p_{\nu} \\ q_{1},\ldots, q_{\mu}}}\frac{a(p_{1})\ldots
a(p_{\nu})\overline{a}(q_{1})\ldots
\overline{a}(q_{\mu})}{p_{1}\ldots q_{\mu}}\biggl(\frac{q_{1}\ldots
q_{\mu}}{p_{1}\ldots p_{\nu}}\biggr)^{it_{n}},
\]
\emph{satisfies the bound} $|S|<(\delta y^{3/2})^{2k}\ln N$.

\vspace{0.2cm}

For a proof, see \cite{Korolev_2010a}.

\vspace{0.2cm}

\textsc{Lemma 3}. \emph{Let} $k\ge 1$ \emph{be an integer},
$y>e^{3}$, \emph{and let} $p_{1},\ldots, p_{k}, q_{1},\ldots, q_{k}$
\emph{take values of prime numbers in the interval} $(1,y]$.
\emph{Then the following relation holds}:
\[
\sum\limits_{p_{1}\ldots p_{k} = q_{1}\ldots
q_{k}}\frac{a(p_{1})\ldots a(p_{k})\overline{a}(q_{1})\ldots
\overline{a}(q_{k})}{\sqrt{p_{1}\ldots
q_{k}}}\,=\,k!\bigl(\sigma_{1}^{k}\,+\,\theta_{k}k^{2}\sigma_{1}^{k-2}\sigma_{2}\bigr),
\]
\emph{where}
\[
\sigma_{j}\,=\,\sum\limits_{p\le
y}\biggl(\frac{|a(p)|^{2}}{p}\biggr)^{\!j},\quad j = 1,2,\quad -1\le
\theta_{k}\le 0,
\]
\emph{and} $\theta_{1}=0$.

\vspace{0.2cm}

For a proof, see \cite{Ghosh_1981},\cite{Karatsuba_Korolev_2006}.

\vspace{0.2cm}

Suppose $x = t_{N}^{0.1\vep}$. For positive $t$ and $y$ we define
\[
V(t) = V_{y}(t) = \frac{1}{\pi}\sum\limits_{p<y}\frac{\sin{(t\ln
p)}}{\sqrt{p}},\quad R(t) = S(t)+V(t).
\]

\vspace{0.2cm}

\textsc{Lemma 4}. \emph{Suppose} $k$ \emph{is an integer such that}
$1\le k\le \tfrac{1}{192}\ln x$, $y = x^{1/(4k)}$, $V(t) =
V_{y}(t)$. \emph{Then}
\[
\sum\limits_{N<n\le N+M}R^{2k}(t_{n}+0)\,\le\,(Ae^{-4}k)^{2k}M,
\]
\emph{where} $A = e^{21}\vep^{-1.5}$.

\vspace{0.2cm}

\textsc{Lemma 5}. \emph{Suppose} $k$ \emph{is an integer}, $1\le
k\le \sqrt{L}$. \emph{Then the following relations hold}:
\begin{align*}
& \sum\limits_{N<n\le
N+M}S^{2k}(t_{n}+0)\,=\,\frac{(2k)!}{k!}\,\frac{ML^{k}}{(2\pi)^{2k}}\,\bigl(1\,+\,\theta A^{k}L^{-0.5}\bigr),\\
& \biggl|\;\sum\limits_{N<n\le
N+M}S^{2k-1}(t_{n}+0)\biggr|\,\le\,\frac{3.5}{\sqrt{B}}\,(Bk)^{k}ML^{k-1},\\
& \sum\limits_{N<n\le
N+M}S^{2k}(t_{n}+0)\,\le\,2\Bigl(\frac{kA}{\pi^{2}e}\Bigr)^{k}ML^{k},
\end{align*}
\emph{where} $B = A^{2}e^{-8}$, \emph{and} $A$ \emph{is defined in
Lemma 4.}

\vspace{0.2cm}

For the proofs of these two lemmas, see \cite{Korolev_2010a} (the
substitution of $t_{n}$ to $t_{n}+0$ does not affect the truth of
the result; the reason is that the above substitution does not
affect to the functions that approximate $S(t)$ in the proofs of
Lemmas 4 and 5).

\vspace{0.2cm}

\textsc{Lemma 6}. \emph{Suppose} $m$ \emph{is an integer}, $1\le m
\le M$. \emph{Then the equality}
\[
t_{n+m}\,-\,t_{n}\,=\,\frac{\pi m}{\vth'(t_{N})}\,+\,\frac{3\theta
M}{N(\ln N)^{2}}
\]
\emph{holds true for} $N<n\le N+M$.

\vspace{0.2cm}

\textsc{Proof.} By Lagrange's mean value theorem, we have
\[
\pi
m\,=\,\vth(t_{n+m})-\vth(t_{n})\,=\,(t_{n+m}\,-\,t_{n})\vth'(\xi),\quad
t_{n+m}\,-\,t_{n}\,=\,\frac{\pi m}{\vth'(\xi)}
\]
for some $\xi$, $t_{n}<\xi <t_{n+m}$. Since $\vth'(t)$, $\vth''(t)$
are monotonic for $t>7$, by the inequality
\[
\frac{t_{N}}{2\pi}\ln{\frac{t_{N}}{2\pi}}\,>\,N
\]
we get:
\begin{multline*}
0\,<\,\frac{1}{\vth'(t_{N})}\,-\,\frac{1}{\vth'(\xi)}\,=\,\frac{\vth'(\xi)\,-\,\vth'(t_{N})}{\vth'(t_{N})\vth'(\xi)}\,\le
\\
\le\,\frac{(t_{N+M}\,-\,t_{N})\vth''(t_{N})}{\bigl(\vth'(t_{N})\bigr)^{2}}\,=\,\frac{\pi
M}{\vth'(\zeta)}\,\frac{\vth''(t_{N})}{\bigl(\vth'(t_{N})\bigr)^{2}}\,<\,\frac{\pi
M\vth''(t_{N})}{\bigl(\vth'(t_{N})\bigr)^{3}}\,<\,\frac{3M}{N(\ln
N)^{2}}.
\end{multline*}
This proves the lemma.

\vspace{0.2cm}

\textsc{Lemma 7}. \emph{Suppose} $0<h_{0}<\tfrac{1}{2}$ \emph{is a
sufficientlty small constant}, $0<h<h_{0}$, $h\ln x > 2$, \emph{and
let}
\[
V(x;h)\,=\,\sum\limits_{p\le x}\frac{\sin^{2}\bigl(\tfrac{1}{2}h\ln
p\bigr)}{p}.
\]
\emph{Then the following relation holds}:
\[
V(x;h)\,=\,\frac{1}{2}\ln{(h\ln x)}\,+\,1.05\theta.
\]

\vspace{0.2cm}

\textsc{Proof}. Setting $y = e^{\lambda/h}$ for some $1<\lambda<2$,
we obtain:
\begin{multline*}
V(x;h)\,=\,\biggl(\;\sum\limits_{p\le y}\,+\,\sum\limits_{y<p\le
x}\;\biggr)\frac{\sin^{2}\bigl(\tfrac{1}{2}h\ln p\bigr)}{p}\,=\\
=\, \sum\limits_{p\le y}\frac{\sin^{2}\bigl(\tfrac{1}{2}h\ln
p\bigr)}{p}\,+\,\frac{1}{2}\sum\limits_{y<p\le x}\frac{1-\cos{(h\ln
p)}}{p}\,=\,V_{1}+\frac{1}{2}(V_{2}-V_{3}).
\end{multline*}
The application of Lemma 1 yields:
\begin{align*}
& V_{1}\,\le\, \Bigl(\frac{h}{2}\Bigr)^{2}\sum\limits_{p\le
y}\frac{(\ln p)^{2}}{p}\,\le\, \Bigl(\frac{h}{2}\Bigr)^{2}(\ln
y)\sum\limits_{p\le y}\frac{\ln p}{p}\,<\,\Bigl(\frac{h}{2}\ln
y\Bigr)^{2} =
\frac{\lambda^{2}}{4},\\
& V_{2} = \biggl(\;\sum\limits_{p\le x}\,-\,\sum\limits_{p\le
y}\;\biggr)\frac{1}{p}\,=\,\ln(h\ln
x)\,-\,\ln\lambda\,+\,2\theta\Bigl(\frac{h}{\lambda}\Bigr)^{2}.
\end{align*}
Now we divide the domain of $p$ in $V_{3}$ into the intervals of the
form $a<p\le b$, where $b\le 2a$, $a = 2^{k}y$, $k = 0,1,2,\ldots$.
Thus we get:
\[
V_{3}\,=\,\RRe\sum\limits_{y<p\le
x}\frac{p^{ih}}{p}\,=\,\RRe\sum\limits_{a}V_{3}(a),\quad
V_{3}(a)\,=\,\sum\limits_{a<p\le b}\frac{p^{ih}}{p}.
\]
Setting
\[
\mathbb{C}(u)\,=\,\sum\limits_{a<p\le
u}\frac{1}{p}\,=\,\ln\ln{u}\,-\,\ln\ln{a}\,+\,\frac{2\theta}{\ln^{2}a}
\]
and applying the Abel's summation formula, we have
\[
V_{3}(a)\,=\,\mathbb{C}(b)b^{ih}\,-\,\int_{a}^{b}\mathbb{C}(u)du^{ih}\,=\,\int_{a}^{b}\frac{u^{ih}du}{u\ln
u}\,+\,\frac{3\theta_{1}}{\ln^{2}a},
\]
and therefore
\[
V_{3}\,=\,\RRe(j)\,+\,3\theta_{2}\sum\limits_{k\ge 0}\frac{1}{(k\ln
2+\ln y)^{2}},\quad j\,=\,\int_{y}^{x}\frac{u^{ih}du}{u\ln u}.
\]
Integration by parts yields:
\[
j\,=\,\frac{1}{ih}\biggl(\frac{x^{ih}}{\ln x}\,-\,\frac{y^{ih}}{\ln
y}\,-\,\int_{y}^{x}u^{ih}d\frac{1}{\ln u}\biggr),\quad
|j|\,\le\,\frac{2}{h\ln y}\,=\,\frac{2}{\lambda}.
\]
Finally, we obtain
\begin{align*}
& |V_{3}|\,\le\,\frac{2}{\lambda}\,+\,\frac{9}{\ln
y}\,=\,\frac{2+9h}{\lambda},\\
& V(x;h)\,=\,\frac{1}{2}\ln{(h\ln x)}\,+\,v(x;h),
\end{align*}
where
\[
|v(x;h)|\,\le\,\frac{\lambda^{2}}{4}\,+\,\frac{1}{\lambda}\,+\,\frac{1}{2}\ln\lambda\,+\,\frac{9h}{\lambda}\,+\,\frac{h^{2}}{\lambda^{2}}.
\]
Setting $\lambda = 1.5$, we arrive at the assertion of the lemma.

\vspace{0.2cm}

\textsc{Lemma 8}. $N_{0}(t)>(0.4+7\cdot 10^{-3})N(t)$ \emph{for}
$t>t_{0}>1$.

\vspace{0.2cm}

For a proof, see \cite{Conrey_1989}.

\pagebreak

\begin{center}
\textbf{\S 3. On mean values of the quantities
$\boldsymbol{S(t_{n+m}+0)\,-\,S(t_{n}+0)}$}
\end{center}

Suppose $m\ge 1$ is an integer. Let us consider the union of $m$
adjacent Gram's intervals $G_{n+1}$, $G_{n+2}, \ldots, G_{n+m}$,
that is the interval $(t_{n}, t_{n+m}]$. By Riemann - von Mangoldt's
formula (\ref{Lab01}), the number of ordinates in last interval
equals to
\begin{multline}
N(t_{n+m}+0)\,-\,N(t_{n}+0)\,=\,\frac{1}{\pi}\bigl(\vth(t_{n+m})\,-\,\vth(t_{n})\bigr)\,+\,S(t_{n+m}+0)\,-\,S(t_{n}+0)\,=\\
=\,m\,+\,S(t_{n+m}+0)\,-\,S(t_{n}+0).\label{Lab11}
\end{multline}
The number of $\gamma_{n}$ that do not exceed a given bound equals
asymptotically to the number of Gram points in the same domain. Then
it is natural to call the number $m$ as an `expected' number of
ordinates of zeros of $\zeta(s)$ in the interval $(t_{n}, t_{n+m}]$.
Hence, the difference
\begin{equation}
S(t_{n+m}+0)\,-\,S(t_{n}+0) \label{Lab12}
\end{equation}
is a deviation of `true' number of ordinates from the `expected'
number.

The below Theorem 1 shows that this deviation often takes a very
large values (of order $\sqrt{\ln m}$, for example). This fact was
observed firstly by Fujii \cite{Fujii_1987} for the case when the
interval of summation is long ($1\le n $ $\le N$) and when $m$
growths with $N$. He proved that the distribution function for the
normalized differences (\ref{Lab12}) tends to Gaussian distribution
as $N\to +\infty$.

\vspace{0.2cm}

\textsc{Theorem 1}. \emph{Let} $k$ \emph{and} $m$ \emph{be an
integers that satisfy the conditions}
\[
k\ge 1, \quad k\vep^{-1}\exp{(\lambda k^{2})}\,\le\,m\,\le\,c\ln N,
\]
\emph{where} $\lambda = (2Be\pi^{2})^{2}$, $B$ \emph{is defined in
Lemma 5}, \emph{and} $c$ \emph{is a sufficiently small absolute
constant. Then the following relation holds}:
\[
\sum\limits_{N<n\le
N+M}\bigl(S(t_{n+m}+0)\,-\,S(t_{n}+0)\bigr)^{2k}\,=\,\frac{(2k)!}{k!}M\Bigl(\frac{1}{2\pi^{2}}\ln\frac{m\vep}{k}\Bigr)^{k}
\biggl(1\,+\,\frac{6\theta
\sqrt{B}\,4^{k}k\sqrt{k}}{\sqrt{\ln{(m\vep\! /k)}}}\biggr).
\]

\vspace{0.2cm}

\textsc{Proof}. Let $x = t_{N}^{0.1\vep}$, $y = x^{1/(4k)}$, $V(t) =
V_{y}(t)$. By lemma 6, $t_{n+m}-t_{n} = h+\vep_{n}$, where
\[
h\,=\,\frac{\pi m}{\vth'(t_{N})},\quad
|\vep_{n}|\,\le\,\frac{3M}{N(\ln N)^{2}}.
\]
By Lagrange's mean value theorem and the inequalities
\[
|V'(t)|\,\le\,\frac{1}{\pi}\sum\limits_{p<y}\frac{\ln
p}{\sqrt{p}}\,<\,\sqrt{y}\,=\,x^{1/(8k)}\,<\,N^{\vep/80}
\]
we get
\[
V(t_{n+m})\,=\,V(t_{n}+h+\vep_{n})\,=\,V(t_{n}+h)\,+\,\vep_{n}V'(t_{n}+h+\theta
\vep_{n})\,=\,V(t_{n}+h)\,+\,\theta_{1}N^{-2/3}.
\]
Hence
\begin{align*}
&
V(t_{n+m})\,-\,V(t_{n})\,=\,\frac{2}{\pi}\,W(t_{n})\,+\,\theta_{2}N^{-2/3},\\
& W(t)\,=\,\frac{1}{2}\sum\limits_{p<y}\frac{\sin{((t+h)\ln
p)}\,-\,\sin{(t\ln
p)}}{\sqrt{p}}\,=\,\sum\limits_{p<y}\frac{\sin{\bigl(\tfrac{1}{2}h\ln
p\bigr)}}{\sqrt{p}}\,\cos{\bigl((t+\tfrac{1}{2}h)\ln p\bigr)}.
\end{align*}
Using the trivial bound $|W(t)|<\sqrt{y}$ and Lagrange's mean value
theorem, we obtain
\begin{align*}
&
\bigl(V(t_{n+m})\,-\,V(t_{n})\bigr)^{2k}\,=\,\biggl(\frac{2}{\pi}\biggr)^{2k}W^{2k}(t_{n})\,+\,\theta
k2^{2k-1}\bigl(|W(t_{n})|^{2k-1}N^{-2/3}\,+\,N^{-4k/3}\bigr)\,=\\
&
=\,\biggl(\frac{2}{\pi}\biggr)^{2k}W^{2k}(t_{n})\,+\,\theta_{1}xN^{-2/3}.
\end{align*}
By summing over $n$ and denoting the corresponding sum by $W_{1}$,
we have
\[
W_{1}\,=\,\sum\limits_{N<n\le
N+M}\bigl(V(t_{n+m})\,-\,V(t_{n})\bigr)^{2k}\,=\,\biggl(\frac{2}{\pi}\biggr)^{2k}W_{2}\,+\,\theta_{2}N^{-1/3},
\]
where
\[
W_{2}\,=\,\sum\limits_{N<n\le N+M}W^{2k}(t_{n}).
\]
Next, we write
$W(t)=\tfrac{1}{2}\bigl(U(t)\,+\,\overline{U}(t)\bigr)$, where
\[
U(t)\,=\,\sum\limits_{p<y}\frac{a(p)}{\sqrt{p}}\,p^{it},\quad
a(p)\,=\,p^{ih/2}\sin{\bigl(\tfrac{1}{2}h\ln p\bigr)}.
\]
Then
\begin{align*}
& W_{2}\,=\,2^{-2k}\sum\limits_{\nu =
0}^{2k}\binom{2k}{\nu}w_{\nu},\\
& w_{\nu}\,=\,\sum\limits_{N<n\le
N+M}\sum\limits_{\substack{p_{1},\ldots, p_{\nu}<y \\
q_{1},\ldots, q_{\mu} < y}}\frac{a(p_{1})\ldots
a(p_{\nu})\overline{a}(q_{1})\ldots
\overline{a}(q_{\mu})}{\sqrt{p_{1}\ldots
q_{\mu}}}\biggl(\frac{p_{1}\ldots p_{\nu}}{q_{1}\ldots
q_{\mu}}\biggr)^{it_{n}},
\end{align*}
where $\mu = 2k-\nu$. By setting $\kappa = \tfrac{1}{4}$, $c = (k\ln
y)(\ln M)^{-1}$, $\delta = 1$ in Lemma 2, we obviously have
$0<c<\tfrac{1}{4} = \tfrac{1}{2}-\kappa$, $N\ge e^{36} =
\exp{(9\kappa^{-1})}$, $y = x^{1/(4k)}>e^{4}$,
\[
\exp{\bigl(3kc^{-1}\bigr)}\,=\,\exp{\Bigl(\frac{3\ln M}{\ln
y}\Bigr)}\,\le\,\exp{\bigl(\tfrac{3}{4}\ln M\bigr)}\,<\,M.
\]
Thus, the conditions of Lemma 2 are satisfied. Hence, for $\nu\ne k$
we have
\[
|w_{\nu}|\,\le\,(y^{3/2})^{2k}\ln N\,=\,x^{3/4}\ln N\,<\,x.
\]
The contribution of the terms of $w_{k}$ that obey the condition
$p_{1}\ldots p_{k}\ne q_{1}\ldots q_{k}$ is estimated as above.
Therefore,
\[
W_{2}\,=\,2^{-2k}\binom{2k}{k}Mw\,+\,\theta 2^{-2k}\sum\limits_{\nu
= 0}^{2k}\binom{2k}{\nu}x\,=\,2^{-2k}\binom{2k}{k}Mw\,+\,\theta x,
\]
where
\[
w\,=\,\sum\limits_{p_{1}\ldots p_{k} = q_{1}\ldots
q_{k}}\frac{a(p_{1})\ldots \overline{a}(q_{k})}{\sqrt{p_{1}\ldots
q_{k}}}.
\]
By Lemma 3, $w =
k!(\sigma_{1}^{k}\,+\,\theta_{k}k^{2}\sigma_{1}^{k-2}\sigma_{2})$,
where
\[
\sigma_{1}\,=\,\sum\limits_{p <
y}\frac{\sin^{2}{\bigl(\tfrac{1}{2}h\ln p\bigr)}}{p},\quad
\sigma_{2}\,\le\,\sum\limits_{p}\frac{1}{p^{2}}\,<\,\frac{1}{2}
\]
and $-1\le \theta_{k}\le 0$. Since
\begin{align*}
& h\ln{y}\,=\,\frac{2\pi
m\ln{y}}{\ln{t_{N}}}\,\bigl(1\,+\,o(1)\bigr)\,=\,\frac{\pi
m\vep}{20k}\,\bigl(1\,+\,o(1)\bigr)\,>\,2,\\
& h\,\le\,\frac{2\pi m}{\ln N}\,\bigl(1\,+\,o(1)\bigr)\,\le\,2\pi
c\,\bigl(1\,+\,o(1)\bigr),
\end{align*}
the conditions of lemma 7 are satisfied for $h_{0} = 7c$ and for
sufficiently small $c$. Hence,
\[
\sigma_{1}\,=\,\frac{1}{2}\ln{(h\ln{y})}\,+\,1.05\theta\,=\,\frac{1}{2}\ln{\frac{\pi
m\vep}{20k}}\,\bigl(1\,+\,o(1)\bigr)\,+\,1.05\theta\,=\,\frac{1}{2}\ln\frac{m\vep}{k}\,+\,2\theta_{1}.
\]
Further, the inequality $\ln{(m\vep/k)}\ge 100k$ implies the
following bounds for $\sigma_{1}$, $w$, $W_{2}$ and $W_{1}$:
\begin{align*}
&
\sigma_{1}^{k}\le\biggl(\frac{1}{2}\ln\frac{m\vep}{k}\biggr)^{k}\!\biggl(1\,+\,\frac{4}{\ln{(m\vep/k)}}\biggr)^{k}\,\le\,
\biggl(\frac{1}{2}\ln\frac{m\vep}{k}\biggr)^{k}\!\biggl(1\,+\,\frac{1}{25k}\biggr)^{k}\!<1.1\biggl(\frac{1}{2}\ln\frac{m\vep}{k}\biggr)^{k},\\
&
w\,\le\,k!\sigma_{1}^{k}\,<\,1.1k!\biggl(\frac{1}{2}\ln\frac{m\vep}{k}\biggr)^{k},\\
&
W_{2}\,\le\,2^{-2k}\frac{(2k)!}{k!}\,M\cdot 1.1\biggl(\frac{1}{2}\ln\frac{m\vep}{k}\biggr)^{k}\,+\,x,\\
& W_{1}\,\le\,\pi^{-2k}\frac{(2k)!}{k!}\,M\cdot
1.1\biggl(\frac{1}{2}\ln\frac{m\vep}{k}\biggr)^{k}\,+\,x\,+\,N^{-1/3}\,<\,2\biggl(\frac{2k}{\pi^{2}e}\biggr)^{k}
\biggl(\ln\frac{m\vep}{k}\biggr)^{k}M.
\end{align*}
Moreover,
\begin{align*}
&
w\,=\,k!\biggl(\biggl(\frac{1}{2}\ln\frac{m\vep}{k}\,+\,2\theta\biggr)^{k}\,+\,\theta_{k}k^{2}\biggl(\frac{1}{2}\ln\frac{m\vep}{k}\,+
\,2\biggr)^{k-2}\biggr)\,=\\
&
=\,k!\biggl(\biggl(\frac{1}{2}\ln\frac{m\vep}{k}\biggr)^{k}\,+\,2\theta
k\biggl(\frac{1}{2}\ln\frac{m\vep}{k}\,+
\,2\biggr)^{k-1}\,+\,\theta_{k}k^{2}\biggl(\frac{1}{2}\ln\frac{m\vep}{k}\,+
\,2\biggr)^{k-2}\biggr)\,=\\
\end{align*}
\begin{align*}
&
=\,k!\biggl(\biggl(\frac{1}{2}\ln\frac{m\vep}{k}\biggr)^{k}\,+\,4\theta
k\biggl(\frac{1}{2}\ln\frac{m\vep}{k}\,+
\,2\biggr)^{k-1}\biggr),\\
&
W_{1}\,=\,\biggl(\frac{2}{\pi}\biggr)^{\!2k}\!\biggl(2^{-2k}\binom{2k}{k}k!M\biggl(\biggl(\frac{1}{2}\ln\frac{m\vep}{k}\biggr)^{k}\!+4\theta
k\biggl(\frac{1}{2}\ln\frac{m\vep}{k}+
2\biggr)^{\!k-1}\biggr)+\theta x\biggr)+\theta
N^{-1/3}\\
&
=\,\frac{(2k)!}{k!}\frac{M}{\pi^{2k}}\biggl(\biggl(\frac{1}{2}\ln\frac{m\vep}{k}\biggr)^{k}\!+4.1\theta
k\biggl(\frac{1}{2}\ln\frac{m\vep}{k}+ 2\biggr)^{\!k-1}\biggr).
\end{align*}
Denoting by $W_{0}$ the initial sum of Lemma and noting that
\[
S(t_{n+m}+0)\,-\,S(t_{n}+0)\,=\,-\bigl(V(t_{n+m})\,-\,V(t_{n})\bigr)\,+\,\bigl(R(t_{n+m}+0)-R(t_{n}+0)\bigr),
\]
we get
\begin{multline*}
\bigl(S(t_{n+m}+0)\,-\,S(t_{n}+0)\bigr)^{2k}\,=\,\bigl(V(t_{n+m})\,-\,V(t_{n})\bigr)^{2k}\,+\,\\
+\,\theta k
2^{2k-1}\bigl(\bigl(R(t_{n+m}+0)-R(t_{n}+0)\bigr)^{2k}\,+\,\bigl|V(t_{n+m})\,-\,V(t_{n})\bigr|^{2k-1}\bigl|R(t_{n+m}+0)-R(t_{n}+0)\bigr|\bigr),\\
\end{multline*}
\[
W_{0}\,=\,\sum\limits_{N<n\le
N+M}\bigl(S(t_{n+m}+0)\,-\,S(t_{n}+0)\bigr)^{2k}\,=\,W_{1}\,+\,\theta_{1}k
2^{2k-1}(W_{3}\,+\,W_{4}),
\]
where
\begin{align*}
& W_{3}\,=\,\sum\limits_{N<n\le
N+M}\bigl(R(t_{n+m}+0)\,-\,R(t_{n}+0)\bigr)^{2k},\\
& W_{4}\,=\,\sum\limits_{N<n\le
N+M}\bigl|V(t_{n+m})\,-\,V(t_{n})\bigr|^{2k-1}\bigl|R(t_{n+m}+0)-R(t_{n}+0)\bigr|.
\end{align*}
By lemma 4,
\begin{multline*}
W_{3}\,\le\,2^{2k-1}\sum\limits_{N<n\le
N+M}\bigl(R^{2k}(t_{n+m}+0)\,+\,R^{2k}(t_{n}+0)\bigr)\,\le\,2^{2k}\sum\limits_{N<n\le
N+2M}R^{2k}(t_{n}+0)\,\le\\
\le\,2^{2k}\cdot 2M(\sqrt{B}k)^{2k}\,=\,2M(2\sqrt{B}k)^{2k}.
\end{multline*}
Further, combining the above bounds for $W_{1}$ and $W_{3}$ with
H\"{o}lder's inequality, we have
\begin{multline*}
W_{4}\,\le\,W_{1}^{1-1/(2k)}W_{3}^{1/(2k)}\,\le\,2^{1/(2k)}\,2\sqrt{B}k\cdot
2^{1-1/(2k)}\biggl(\frac{2k}{\pi^{2}e}\ln\frac{m\vep}{k}\biggr)^{k-1/2}M\,=\\
=\,4\sqrt{B}k\biggl(\frac{2k}{\pi^{2}e}\ln\frac{m\vep}{k}\biggr)^{k-1/2}M.
\end{multline*}
Therefore,
\begin{multline*}
k 2^{2k-1}\bigl(W_{3}\,+\,W_{4}\bigr)\,\le\,k 2^{2k-1}\cdot
4Mk\sqrt{B}\biggl(\frac{2k}{\pi^{2}e}\ln\frac{m\vep}{k}\,\biggr)^{\!k-1/2}M\,\times\\
\biggl(1\,+\,\biggl(\frac{k\sqrt{\lambda}}{\ln{(m\vep/k)}}\biggr)^{k-1/2}\biggr)\,<\,4.1Mk^{2}\sqrt{B}\biggl(\frac{8k}{\pi^{2}e}\ln\frac{m\vep}{k}\biggr)^{k-1/2}.
\end{multline*}
Finally, we get
\begin{align*}
&
W_{0}=\frac{(2k)!}{k!}\,\frac{M}{\pi^{2k}}\biggl(\biggl(\frac{1}{2}\ln\frac{m\vep}{k}\biggr)^{\!k}\!\!+4.1k\theta_{1}
\biggl(\frac{1}{2}\ln\frac{m\vep}{k}+2\biggr)^{\!k-1}\biggr)\,+\\
&+\,4.1\theta_{2}\sqrt{B}k^{2}\biggl(\frac{8k}{\pi^{2}e}\ln\frac{m\vep}{k}\biggr)^{\!k-1/2}\!\!M\,=\,\frac{(2k)!}{k!}\,M\,\biggl(\frac{1}{2\pi^{2}}\ln\frac{m\vep}{k}\biggr)^{k}\bigl(1\,+\,\theta(\delta_{1}+\delta_{2})\bigr),
\end{align*}
where
\begin{align*}
& \delta_{1}\,=\,4.1k\biggl(1\,+\,\frac{4}{\ln{(m\vep/k)}}\biggr)^{k-1}\frac{2}{\ln{(m\vep/k)}}\,<\,\frac{8.6k}{\ln{(m\vep/k)}}, \\
&
\delta_{2}\,=\,\frac{\pi^{2k}k!}{(2k)!}\,\frac{4.1\sqrt{B}k^{2}}{\bigl(\tfrac{1}{2}\ln{(m\vep/k)}\bigr)^{k}}\biggl(\frac{8k}{\pi^{2}e}\ln\frac{m\vep}{k}\biggr)^{k-1/2}
\,\le\,5.8\sqrt{B}\,\frac{4^{k}k\sqrt{k}}{\sqrt{\ln{(m\vep/k)}}}.
\end{align*}
It remains to note that
\[
\delta_{1}\,+\,\delta_{2}\,<\,6\sqrt{B}\,\frac{4^{k}k\sqrt{k}}{\sqrt{\ln{(m\vep/k)}}}.
\]
Thus the theorem is proved.

\vspace{0.2cm}

\textsc{Corollary.} \emph{Suppose} $k$ \emph{and} $m$ \emph{are
integers such that}
\[
1\le k\le 0.1\ln\ln\ln N, \quad k\vep^{-1}\exp{(\varkappa)}<m\le
c\ln N,
\]
\emph{where} $c$ \emph{is sufficiently small absolute constant},
\emph{and} $\varkappa$ \emph{is a maximum of the numbers}
$(2Be\pi^{2}k)^{2}$ \emph{and} $12^{2}B\,4^{2k}k^{3}$. \emph{Then
the following inequality holds}:
\[
\sum\limits_{N<n\le
N+M}\bigl(S(t_{n+m}+0)\,-\,S(t_{n}+0)\bigr)^{2k}\,>\,\frac{M}{2}\frac{(2k)!}{k!}\biggl(\frac{1}{2\pi^{2}}\ln\frac{m\vep}{k}\biggr)^{k}.
\]
\emph{In particular}, \emph{if}
$\vep^{-1}\exp{\bigl((2Be\pi^{2})^{2}\bigr)}<m\le c\ln N$,
\emph{then}
\[
\sum\limits_{N<n\le
N+M}\bigl(S(t_{n+m}+0)\,-\,S(t_{n}+0)\bigr)^{2}\,>\,\frac{M}{2\pi^{2}}\ln{(m\vep)},
\]
\emph{and if} $m = [\mu]+1$, $\mu =
\vep^{-1}\exp{(e^{76}\vep^{-6})}$, \emph{then}
\[
\sum\limits_{N<n\le
N+M}\bigl(S(t_{n+m}+0)\,-\,S(t_{n}+0)\bigr)^{2}\,>\,1.01e^{73}\vep^{-6}M.
\]

\vspace{0.2cm}

\textsc{Theorem 2}. \emph{Let} $k$ \emph{be an integer such that}
$1\le k\le \tfrac{1}{192}\ln x$. \emph{Then the following
ine\-qu\-a\-li\-ty holds}:
\[
\sum\limits_{N<n\le
N+M}\bigl(S(t_{n}+0)\,-\,S(t_{n-1}+0)\bigr)^{2k}\,\le\,2Mk(4k\sqrt{B})^{2k}.
\]

\vspace{0.2cm}

\textsc{Proof.} Using the same arguments as above, in the case $m =
1$ we obtain
\[
W_{1}\,=\,\sum\limits_{N<n\le
N+M}\bigl(V(t_{n})\,-\,V(t_{n-1})\bigr)^{2k}\,\le\,\pi^{-2k}\,\frac{(2k)!}{k!}\,M\sigma_{1}^{k}\,+\,x,
\]
where
\[
\sigma_{1}\,=\,\sum\limits_{p<y}\frac{\sin^{2}{\bigl(\tfrac{1}{2}h\ln
p\bigr)}}{p},\quad h\,=\,\frac{\pi}{\vth'(t_{N})}.
\]
Using the relation
\[
\frac{h}{2}\ln
y\,=\,\frac{\pi}{\ln{t_{N}}+O(1)}\,\frac{\ln{x}}{4k}\,<\,\frac{\ln{x}}{k\ln{t_{N}}}\,<\,\frac{\vep}{10k},
\]
and applying Lemma 1, we get
\begin{align*}
&
\sigma_{1}\,\le\,\frac{h^{2}}{4}\sum\limits_{p<y}\frac{\ln^{2}p}{p}\,<\,\biggl(\frac{h}{2}\ln
y\biggr)^{2}\,<\,\biggl(\frac{\vep}{10k}\biggr)^{2},\\
& W_{1}\,<\,\frac{(2k)!}{k!}\,M\biggl(\frac{\vep}{10\pi
k}\biggr)^{2k}\,+\,x\,<\,\frac{3}{2}\biggl(\frac{\vep^{2}}{25\pi^{2}ek}\biggr)^{k}M\,+\,x\,<\,\vep^{2k}M.
\end{align*}
The application of H\"{o}lder's inequality to the initial sum
$W_{0}$ of the lemma yields:
\[
W_{0}\,=\,\sum\limits_{N<n\le
N+M}\bigl(S(t_{n}+0)\,-\,S(t_{n-1}+0)\bigr)^{2k}\,\le\,W_{1}\,+\,k2^{2k-1}\bigl(W_{3}\,+\,W_{1}^{1-1/(2k)}W_{3}^{1/(2k)}\bigr),
\]
where the sum
\[
W_{3}\,=\,\sum\limits_{N<n\le
N+M}\bigl(R(t_{n}+0)\,-\,R(t_{n-1}+0)\bigr)^{2k}
\]
was estimated in the proof of previous theorem. Using the above
bounds of $W_{1}$ и $W_{3}$, we obtain:
\begin{align*}
&
k2^{2k-1}\bigl(W_{3}\,+\,W_{1}^{1-1/(2k)}W_{3}^{1/(2k)}\bigr)\,\le\,1.5Mk(4k\sqrt{B})^{2k},\\
& W_{0}\,\le\,
\vep^{2k}M\,+\,1.5Mk(4k\sqrt{B})^{2k}\,<\,2Mk(4k\sqrt{B})^{2k}.
\end{align*}
The theorem is proved.

\vspace{0.2cm}

The above theorems imply the lower bound for the `first moment' of
the differences $S(t_{n}+0)-S(t_{n-1}+0)$.

\vspace{0.2cm}

\textsc{Theorem 3}. \emph{There exists a positive constant} $c_{1} =
c_{1}(\vep)$ \emph{such that}
\[
\sum\limits_{N<n\le
N+M}\bigl|S(t_{n}+0)\,-\,S(t_{n-1}+0)\bigr|\,>\,c_{1}M.
\]

\vspace{0.2cm}

\textsc{Proof}. Let us put for brevity $r(n) =
S(t_{n}+0)-S(t_{n-1}+0)$ and denote by $V_{k}$ the sum of $k$th
powers of $|r(n)|$. Further, let $m = [\mu]+1$, where $\mu =
\vep^{-1}\exp{\bigl(e^{76}\vep^{-6}\bigr)}$. Combining the identity
\[
S(t_{n+m}+0)\,-\,S(t_{n}+0)\,=\,r(n+1)\,+\,r(n+2)\,+\,\ldots\,+\,r(n+m)
\]
with Cauchy's inequality, we get
\[
\bigl(S(t_{n+m}+0)\,-\,S(t_{n}+0)\bigr)^{2}\,\le\,m\sum\limits_{\nu
= 1}^{m}r^{2}(n+\nu).
\]
Summing both parts of the above relation over $n$, we have
\begin{multline*}
\sum\limits_{N<n\le
N+M}\bigl(S(t_{n+m}+0)\,-\,S(t_{n}+0)\bigr)^{2}\,\le\,m\sum\limits_{\nu
= 1}^{m}\sum\limits_{N<n\le N+M}r^{2}(n+\nu)\,\le\\
\le\,m\sum\limits_{\nu = 1}^{m}\sum\limits_{N<n\le
N+M+m}r^{2}(n)\,=\,m^{2}\sum\limits_{N<n\le N+M+m}r^{2}(n).
\end{multline*}
Hence, by the Corollary of Theorem 1, we get:
\[
\sum\limits_{N<n\le N+M+m}r^{2}(n)\,\ge\,m^{-2}\sum\limits_{N<n\le
N+M}\bigl(S(t_{n+m}+0)\,-\,S(t_{n}+0)\bigr)^{2}\,\ge\,1.01e^{73}\vep^{-6}m^{-2}\,M.
\]
Since $|r(n)|\le |S(t_{n})|+|S(t_{n-1})|\le 18\ln N$ for $N<n\le
N+M+m$ (see \cite{Karatsuba_Korolev_2005}), we have:
\[
V_{2}\,=\,\sum\limits_{N<n\le
N+M}r^{2}(n)\,\ge\,1.01e^{73}\vep^{-6}m^{-2}\,M\,-\,m(18\ln
N)^{2}\,>\,c_{2}M,
\]
where $c_{2} = e^{73}\vep^{-6}m^{-2}\,M$. Further, the application
of H\"{o}lder's inequality to the sum $V_{2}$ yields
\begin{multline*}
V_{2}\,=\,\sum\limits_{N<n\le
N+M}|r(n)|^{2/3}|r(n)|^{4/3}\,\le\,\biggl(\sum\limits_{N<n\le
N+M}|r(n)|\biggr)^{2/3}\biggl(\sum\limits_{N<n\le
N+M}|r(n)|^{4}\biggr)^{1/3}\,=\\
=\,V_{1}^{2/3}V_{4}^{1/3}.
\end{multline*}
Therefore, $V_{1}\ge V_{2}^{3/2}V_{4}^{-1/2}$. Using both the above
bound for $V_{2}$ and the inequality of Theorem 2 for $k = 2$, we
obtain:
\begin{align*}
& V_{4}\,\le\,2^{14}B^{2}M\,=\,c_{4}M,\\
& V_{1}\,=\,\sum\limits_{N<n\le
N+M}|r(n)|\,\ge\,\frac{(c_{2}M)^{3/2}}{(c_{4}M)^{1/2}}\,=\,c_{1}M,\quad
c_{1}\,=\,c_{2}^{3/2}c_{4}^{-1/2}.
\end{align*}
Theorem is proved.

\vspace{0.2cm}

The following assertion is an analogue (for the short interval of
summation) of corollary of Selberg's theorem cited in \S 1.

\vspace{0.2cm}

\textsc{Theorem 4}. \emph{There exist positive constants} $K_{1}$
\emph{and} $K_{2}$ \emph{such that for} $N<n\le N+M$, \emph{there
are more than} $K_{1}M$ \emph{cases when the interval} $G_{n}$
\emph{does not contain any ordinate of a zero of} $\zeta(s)$,
\emph{and more than} $K_{2}M$ \emph{cases when the interval} $G_{n}$
\emph{contains at least two ordinates}, \emph{i.e.}
\[
\nu_{0}(N+M)\,-\,\nu_{0}(N)\,\ge\,K_{1}M,\quad \sum\limits_{k\ge
2}\bigl(\nu_{k}(N+M)\,-\,\nu_{k}(N)\bigr)\,\ge\,K_{2}M.
\]

\vspace{0.2cm}

\textsc{Proof}. Since $r(n) = S(t_{n}+0)-S(t_{n-1}+0)$ is an integer
and $r(n)\ge -1$ for any $n$, then the equality (\ref{Lab11})
implies that the interval $G_{n}$ does not contain any ordinate iff
$r(n) = -1$ and contains more than one ordinate iff $r(n)\ge 1$. In
other words, the number $M_{1}$ of `empty' Gram's intervals is equal
to the number of $n$ such that $r(n)$ is negative, and the number
$M_{2}$ of Gram's intervals that contain two or more ordinates is
equal to the number of positive $r(n)$.

Using the relation
\begin{equation*}
\frac{1}{2}\bigl(|r(n)|\,-\,r(n)\bigr)\,=\,
\begin{cases}
1, & \text{if} \;\; r(n) < 0,\\
0, & \text{if} \;\; r(n) > 0,
\end{cases}
\end{equation*}
by the bound of Theorem 3 we get
\begin{multline*}
M_{1}\,=\,\sum\limits_{N<n\le
N+M}\frac{1}{2}\bigl(|r(n)|\,-\,r(n)\bigr)\,=\,\frac{1}{2}\sum\limits_{N<n\le
N+M}\bigl|S(t_{n}+0)\,-\,S(t_{n-1}+0)\bigr|\,-\\
-\,\frac{1}{2}\sum\limits_{N<n\le
N+M}\bigl(S(t_{n}+0)\,-\,S(t_{n-1}+0)\bigr)\,\ge\,\frac{c_{1}}{2}M\,-\,9\ln
N\,>\,K_{1}M,
\end{multline*}
where $K_{1} = \tfrac{2}{5}c_{1}$. Further, $M_{2}$ is equal to the
number of non-zero terms of the sum
\[
W\,=\,\sum\limits_{N<n\le
N+M}\frac{1}{2}\bigl(|r(n)|\,+\,r(n)\bigr).
\]
The application of Theorem 2 and Cauchy's inequality yields:
\[
W\,\le\,\sqrt{M_{2}}\sqrt{\sum\limits_{N<n\le
N+M}r^{2}(n)}\,\le\,\sqrt{M_{2}}\sqrt{32BM}.
\]
Since $W>K_{1}M$, then $M_{2}\ge K_{2}M$, where
\[
K_{2}\,=\,\frac{K_{1}^{2}}{32B}.
\]
Theorem is proved.

\textsc{Remark}. The constants $K_{1}$, $K_{2}$ defined above are
too small. It's easy to see that they do not exceed
$\exp{\bigl(-e^{75}\vep^{-6}\bigr)}$. At the same time, the
calculations of zeros of $\zeta(s)$ shows that probably $K_{1}>0.1$,
$K_{2}>0.1$. Thus, it is of some interest to prove the analogue of
Theorem 4 with the constants $K_{1}$ and $K_{2}$ of order $0.001 -
0.01$.

\begin{center}
\textbf{\S 4. The alternating sums connected with the function
$\boldsymbol{S(t)}$.}
\end{center}

Here we study the sums of the following type:
\[
T_{k}\,=\,\sum\limits_{N<n\le
N+M}S^{k}(t_{n}+0)\bigl(S(t_{n}+0)\,-\,S(t_{n-1}+0)\bigr).
\]
Theorem 4 implies that the difference $r(n) =
S(t_{n}+0)-S(t_{n-1}+0)$ is negative for positive proportion of $n$,
$N<n\le N+M$. At the same time, this difference is positive for
positive proportion of $n$. Hence, the sums $T_{k}$ are alternating
(at least, for even $k$).

The direct application of Cauchy's inequality does not allow us to
take into account the oscillation in the sum $T_{k}$, and therefore
does not allow us to obtain non-trivial bound for $T_{k}$. Thus, for
example,
\[
|T_{2k}|\le\Bigl(\sum\limits_{N<n\le
N+M}S^{4k}(t_{n}+0)\Bigr)^{1/2}\Bigl(\sum\limits_{N<n\le
N+M}r^{2}(n)\Bigr)^{1/2}\,\ll\,_{k}\,\sqrt{ML^{2k}}\,\sqrt{M}\,\ll\,_{k}\,ML^{k},
\]
and this bound is trivial. Hence, we need to use some additional
arguments.

\vspace{0.2cm}

\textsc{Theorem 5}. \emph{Suppose} $k$ \emph{is an integer such
that} $1\le k\le \sqrt{L}$. \emph{Then the following estimations
hold}:
\[
|T_{2k-1}|\,<\,0.02(Ak)^{k+1}ML^{k-1},\quad\quad
|T_{2k}|\,<\,0.02(10A)^{k+1}\,\frac{(2k)!}{k!}\frac{ML^{k-1/2}}{(2\pi)^{2k}}.
\]

\vspace{0.2cm}

\textsc{Proof.} Consider first the sum $T_{2k-1}$. Setting $a =
S(t_{n}+0)$, $b = r(n) = S(t_{n}+0)-S(t_{n-1}+0)$ in the
easy-to-check identity
\[
(a-b)^{2k}\,=\,a^{2k}\,-\,2ka^{2k-1}b\,+\,\theta
k^{2}2^{2k-2}\bigl(a^{2k-2}b^{2}\,+\,b^{2k}\bigr),
\]
after some obvious transformations we get:
\[
2k\,
S^{2k-1}(t_{n}+0)r(n)\,=\,S^{2k}(t_{n}+0)\,-\,S^{2k}(t_{n-1}+0)\,+\,\theta
k^{2}2^{2k-2}\bigl(r^{2k}(n)\,+\,S^{2k-2}(t_{n}+0)r^{2}(n)\bigr).
\]
Summing over $n$, we obtain
\[
2k\,T_{2k-1}\,=\,S^{2k}(t_{N+M})\,-\,S^{2k}(t_{N})\,+\,\theta\,k^{2}2^{2k-2}\bigl(W_{1}\,+\,W_{2}\bigr),
\]
where
\[
W_{1}\,=\,\sum\limits_{N<n\le N+M}r^{2k}(n),\quad
W_{2}\,=\,\sum\limits_{N<n\le N+M}S^{2k-2}(t_{n}+0)r^{2}(n).
\]
By theorem 2,
\[
W_{1}\,\le\,2Mk(4k\sqrt{B})^{2k}.
\]
Using both the estimation of Lemma 5 and H\"{o}lder's inequality, we
have:
\[
W_{2}\,\le\,\biggl(\;\sum\limits_{N<n\le
N+M}S^{2k}(t_{n}+0)\biggr)^{1-1/k}W_{1}^{1/k}\,\le\,2\sqrt[3\;]{3}(4k\sqrt{B})^{2}\biggl(\frac{kAL}{\pi^{2}e}\biggr)^{k-1}M.
\]
Hence
\begin{align*}
&
k2^{2k-2}(W_{1}+W_{2})\,\le\,k2^{2k-2}\,2\sqrt[3\;]{3}(4k\sqrt{B})^{2}\biggl(\frac{kAL}{\pi^{2}e}\biggr)^{k-1}M\,\biggl(1\,+\,\frac{k}{\sqrt[3\;]{3}}
\biggl(\frac{kA}{4L}\biggr)^{k-1}\biggr)\,<\\
& <\,2\cdot
4^{2}\sqrt[3\;]{3}e^{-8}k^{3}A^{2}\biggl(\frac{4kAL}{\pi^{2}e}\biggr)^{k-1}\,<\,\frac{1}{30}k^{3}A^{2}(0.15kAL)^{k-1}M\,\le\,\frac{1}{30}(Ak)^{k+1}ML^{k-1},
\end{align*}
and therefore
\[
|T_{2k-1}|\,\le\,\frac{1}{60}(Ak)^{k+1}ML^{k-1}\,+\,\frac{1}{k}(9\ln
N)^{2k}\,<\,\frac{1}{50}(Ak)^{k+1}ML^{k-1}.
\]

Now we consider the sum $T_{2k}$. Setting $a = S(t_{n}+0)$, $b =
r(n)$ in the identity
\[
(a-b)^{2k+1}\,=\,a^{2k+1}\,-\,(2k+1)a^{2k}b\,+\,\theta
k(2k+1)2^{2k-2}\bigl(|b|^{2k+1}\,+\,|a|^{2k-1}b^{2}\bigr),
\]
after some transformations we get:
\[
(2k+1)T_{2k+1}\,=\,S^{2k+1}(t_{N+M})\,-\,S^{2k+1}(t_{N})\,+\,\theta
k(2k+1)2^{2k-2}\bigl(W_{1}\,+\,W_{2}\bigr),
\]
where
\[
W_{1}\,=\,\sum\limits_{N<n\le N+M}|r(n)|^{2k+1},\quad
W_{2}\,=\,\sum\limits_{N<n\le N+M}|S(t_{n}+0)|^{2k-1}r^{2}(n).
\]
We have
\begin{align*}
& W_{1}\,=\,\sum\limits_{N<n\le
N+M}|r(n)|^{2k-1}r^{2}(n)\,\le\,\biggl(\;\sum\limits_{N<n\le
N+M}r^{2k}(n)\biggr)^{\!1-1/(2k)}\!\biggl(\;\sum\limits_{N<n\le
N+M}r^{4k}(n)\biggr)^{\!1/(2k)}\\
&
\le\,\bigl(2k(4k\sqrt{B})^{2k}\bigr)^{1-1/(2k)}\,\bigl(4k(8k\sqrt{B})^{4k}\bigr)^{1/(2k)}M\,<\,8\sqrt{2}kM\bigl(4k\sqrt{B}\bigr)^{2k+1}\,=\\
& =\,\frac{(2k)!}{k!}\frac{M}{(2\pi)^{2k}}\,(AL)^{k-1/2}\delta_{1},
\end{align*}
where
\begin{align*}
&
\delta_{1}\,=\,\frac{(2\pi)^{2k}k!}{(2k)!}\,\frac{8\sqrt{2}k(4k\sqrt{B})^{2k+1}}{(AL)^{k-1/2}}\,\le\,
k^{2}A\sqrt{AL}\biggl(\frac{(4\pi)^{2}Ak}{e^{7}L}\biggr)^{k}\,<\,1.
\end{align*}

Next, the application of Lemmas 4 and 5 yields:
\begin{align*}
& W_{2}\,\le\,\biggl(\;\sum\limits_{N<n\le
N+M}S^{2k}(t_{n}+0)\biggr)^{1-1/(2k)}\biggl(\;\sum\limits_{N<n\le
N+M}r^{4k}(n)\biggr)^{1/(2k)}\,\le\,\\
&
\le\,\biggl(\frac{1.1}{A}\,\frac{(2k)!}{(2\pi)^{2k}k!}\,(AL)^{k}\biggr)^{1-1/(2k)}\!(4k)^{1/(2k)}\bigl(8k\sqrt{B}\bigr)^{2}\,\le\,
\frac{(2k)!}{k!}\frac{M}{(2\pi)^{2k}}(AL)^{k-1/2}\delta_{2},
\end{align*}
where
\[
\delta_{2}\,=\,2\pi\biggl(\frac{k!}{(2k)!}\biggr)^{1/(2k)}\!\sqrt{\frac{1.1}{A}}\,(4k)^{1/(2k)}\bigl(8k\sqrt{B}\bigr)^{2}\,<\,0.01(kA)^{3/2}.
\]
Let us note that
\begin{align*}
& k 2^{2k-2}\bigl(W_{1}\,+\,W_{2}\bigr)\,\le\,k
2^{2k-2}\,\frac{(2k)!}{k!}\frac{M}{(2\pi)^{2k}}\,(AL)^{k-1/2}(\delta_{1}+\delta_{2})\,<\\
& <\,\frac{(2k)!}{k!}\frac{ML^{k-1/2}}{(2\pi)^{2k}}\cdot
0.01(10A)^{k+1}.
\end{align*}
Therefore,
\begin{align*}
& (2k+1)|T_{2k}|\,\le\,(2k+1)\cdot
0.01(10A)^{k+1}\,\frac{(2k)!}{k!}\frac{ML^{k-1/2}}{(2\pi)^{2k}}\,+\,2(9\ln
N)^{2k+1},\\
&
|T_{2k}|\,<\,0.02(10A)^{k+1}\frac{(2k)!}{k!}\frac{ML^{k-1/2}}{(2\pi)^{2k}}.
\end{align*}
Theorem is proved.

\vspace{0.2cm}

The following lemma is necessary for the proof of Selberg's
formulas.

\vspace{0.2cm}

\textsc{Lemma 9.} \emph{Let} $k$ \emph{and} $n$ \emph{be an
arbitrary natural numbers and suppose that the interval} $G_{n} =
(t_{n-1},t_{n}]$ \emph{contains} $(r+1)$ \emph{ordinates of zeros
of} $\zeta(s)$, $r = r(n)\ge -1$. \emph{Then the following relations
hold}:
\begin{align}
& \sum\limits_{t_{n-1}<\gamma_{m}\le
t_{n}}\Delta_{m}^{2k}\,=\,(r+1)\Delta^{2k}(n)\,+\,\theta_{1}k
2^{2k}\bigl(|\Delta(n)|^{2k-1}r^{2}\,+\,|r|^{2k+1}\bigr),\label{Lab13}\\
& \sum\limits_{t_{n-1}<\gamma_{m}\le
t_{n}}\Delta_{m}^{2k-1}\,=\,-(r+1)\Delta^{2k-1}(n)\,+\,\theta_{2}k
2^{2k}\bigl(\Delta^{2k-2}(n)r^{2}\,+\,r^{2k}\bigr).\label{Lab14}
\end{align}

\vspace{0.2cm}

\textsc{Proof}. First we consider the case $r = -1$. Then $G_{n}$
does not contain any ordinate, and the sums in the left-hand sides
of (\ref{Lab13}),(\ref{Lab14}) are empty. Thus the assertion of
lemma is true for $\theta_{1} = \theta_{2} = 0$.

Now let us consider the case $r\ge 0$. Suppose that the inequalities
\begin{equation*}
\gamma_{s-1}\,\le\,t_{n-1}\,<\,\gamma_{s}\,\le\,\gamma_{s+1}\,\le\,\ldots\,\le\,\gamma_{s+r}\,\le\,t_{n}\,<\,\gamma_{s+r+1}
\end{equation*}
hold for some $s\ge 1$. Then $\Delta(n) = S(t_{n}+0) =
N(t_{n}+0)-\pi^{-1}\vth(t_{n})-1 = s+r-n$, and hence
\begin{align*}
& \Delta_{s}\,=\,n-s\,=\,r-\Delta(n),\\
& \Delta_{s+1}\,=\,n-s-1\,=\,r-1-\Delta(n),\\
& \ldots \\
& \Delta_{s+r}\,=\,n-s-r\,=\,-\Delta(n).
\end{align*}
Therefore,
\begin{multline*}
\sum\limits_{t_{n-1}<\gamma_{m}\le
t_{n}}\Delta_{m}^{2k}\,=\,\sum\limits_{m =
s}^{s+r}\Delta_{m}^{2k}\,=\\
=\,\sum\limits_{j =
0}^{r}\bigl(j\,-\,\Delta(n)\bigr)^{2k}\,=\,\sum\limits_{j =
0}^{r}\bigl(\Delta^{2k}(n)\,+\,\theta\,2kj\,2^{2k-2}(|\Delta(n)|^{2k-1}\,+\,j^{2k-1})\bigr)\,=\\
 =\,(r+1)\Delta^{2k}(n)\,+\,\theta_{1}k
2^{2k}\bigl(|\Delta(n)|^{2k-1}r^{2}\,+\,r^{2k+1}\bigr).
\end{multline*}
The proof of (\ref{Lab14}) follows the same arguments. Lemma is
proved.

\vspace{0.2cm}

The below theorem reduces the calculation of the sums
(\ref{Lab08}),(\ref{Lab09}) to the calculation of the sums of the
quantities $\Delta^{k}(n) = S^{k}(t_{n}+0)$ (see Lemma 5).

\vspace{0.2cm}

\textsc{Theorem 6}. \emph{Let} $k$ \emph{be an integer such that}
$1\le k\le \sqrt{L}$. \emph{Then the following relations hold}:
\begin{align*}
& \sum\limits_{N<n\le
N+M}\Delta_{n}^{2k}\,=\,\frac{(2k)!}{k!}\,\frac{ML^{k}}{(2\pi)^{2k}}\,\bigl(1\,+\,\theta
(10A)^{k+1}L^{-0.5}\bigr),\\
& \biggl|\,\sum\limits_{N<n\le
N+M}\Delta_{n}^{2k-1}\biggr|\,\le\,e^{9}(Bk)^{k}ML^{k-1}.
\end{align*}

\vspace{0.2cm}

\textsc{Proof.} The numbers $\mu$, $\nu$ are uniquely defined by the
conditions
\[
\gamma_{\mu}\,\le\,t_{N}\,<\,\gamma_{\mu+1},\quad
\gamma_{\nu}\,\le\,t_{N+M}\,<\,\gamma_{\nu+1}.
\]
Let us consider the sum
\[
V\,=\,\sum\limits_{\mu<m\le \nu}\Delta_{m}^{2k}.
\]
Using the definition of $\Delta(n)$, we obtain:
\begin{align*}
& \mu\,=\,N(t_{N}+0)\,=\,\pi^{-1}\vth(t_{N})\,+\,1\,+\,S(t_{N}+0) =
N+\Delta(N),\\
& |\mu-N|\,=\,|\Delta(N)|\,<\,9\ln{N},
\end{align*}
and, similarly, $|\nu\,-\,(N+M)|\,=\,|\Delta(N+M)|\,<\,9\ln{N}$.
These inequalities and the bound $|\Delta_{m}|<9\ln{N}$ (see Lemma 7
and a posterior remark in \cite{Korolev_2010b}) imply that the
difference between the sum $V$ and the sum
\[
\sum\limits_{N<m\le N+M}\Delta_{m}^{2k}
\]
does not exceed in modulus
\[
(18\ln N + 1)(9\ln N)^{2k}\,<\,3(9\ln N)^{2k+1}.
\]

On the other hand, Lemma 9 implies that
\begin{align*}
& V\,=\,\sum\limits_{N<n\le N+M}\sum\limits_{t_{n-1}<\gamma_{m}\le
t_{n}}\Delta_{m}^{2k}\,=\\
& =\,\sum\limits_{N<n\le N+M}\bigl((r(n)+1)\Delta^{2k}(n)\,+\,\theta
k\,2^{2k}(|\Delta(n)|^{2k-1}r^{2}(n)\,+\,|r(n)|^{2k+1})\bigr)\,=\\
& =\,\sum\limits_{N<n\le
N+M}\Delta^{2k}(n)\,+\,T_{2k}\,+\,\theta_{1}
k\,2^{2k}\bigl(W_{1}+W_{2}\bigr),
\end{align*}
where $r(n) = S(t_{n}+0)-S(t_{n-1}+0)$ and
\begin{align*}
& W_{1}=\sum\limits_{N<n\le N+M}|r(n)|^{2k+1}=\sum\limits_{N<n\le
N+M}|S(t_{n}+0)-S(t_{n-1}+0)|^{2k+1},\\
& W_{2}=\sum\limits_{N<n\le
N+M}|\Delta(n)|^{2k-1}r^{2}(n)=\sum\limits_{N<n\le
N+M}|S(t_{n}+0)|^{2k-1}\bigl(S(t_{n}+0)-S(t_{n-1}+0)\bigr)^{2}.
\end{align*}
Proving Theorem 5, we found that
\begin{align*}
&
|T_{2k}|\,\le\,0.02(10A)^{k+1}\,\frac{(2k)!}{k!}\,\frac{ML^{k-0.5}}{(2\pi)^{2k}},\\
&
k\,2^{2k}\bigl(W_{1}\,+\,W_{2}\bigr)\,<\,0.08(10A)^{k+1}\,\frac{(2k)!}{k!}\,\frac{ML^{k-0.5}}{(2\pi)^{2k}}.
\end{align*}
Using both these inequalities and the assertion of Lemma 5, we get:
\begin{multline*}
\sum\limits_{N<m\le N+M}\Delta_{m}^{2k}=\sum\limits_{N<n\le
N+M}\Delta^{2k}(n)\,+\,0.1\theta_{1}
(10A)^{k+1}\,\frac{(2k)!}{k!}\,\frac{ML^{k-0.5}}{(2\pi)^{2k}}\,+\,3\theta_{2}(9\ln N)^{2k+1}\,=\\
=\,\frac{(2k)!}{k!}\frac{ML^{k}}{(2\pi)^{2k}}\,\bigl(1\,+\,\theta_{3}A^{k}L^{-0.5}\,+\,0.2\theta_{4}
(10A)^{k+1}L^{-0.5}\bigr)\,=\\
=\,\frac{(2k)!}{k!}\frac{ML^{k}}{(2\pi)^{2k}}\bigl(1\,+\,\theta
(10A)^{k+1}L^{-0.5}\bigr).
\end{multline*}
Applying the same arguments, we obtain
\begin{align*}
& V\,=\,\sum\limits_{\mu<m\le
\nu}\Delta_{m}^{2k-1}\,=\,-\sum\limits_{N<n\le
N+M}\Delta^{2k-1}(n)\,-\,T_{2k-1}\,+\,\theta
k\,2^{2k}\bigl(V_{1}\,+\,V_{2}\bigr),
\end{align*}
where
\[
V_{1}\,=\,\sum\limits_{N<n\le N+M}r^{2k}(n),\quad
\sum\limits_{N<n\le N+M}\Delta^{2k-2}(n)r^{2}(n).
\]
Proving Theorem 5, we also found that
\[
k\,2^{2k}\bigl(V_{1}\,+\,V_{2}\bigr)\,\le\,0.14(Ak)^{k+1}ML^{k-1},\quad
|T_{2k-1}|\,<\,0.02(Ak)^{k+1}ML^{k-1}.
\]
By these estimations and by the inequalities of Lemma 5 we have:
\begin{multline*}
\biggl|\sum\limits_{N<m\le
N+M}\Delta_{m}^{2k-1}\biggr|\,\le\,\biggl|\sum\limits_{N<n\le
N+M}\Delta^{2k-1}(n)\biggr|\,+\,|T_{2k-1}|\,+\,k\,2^{2k}\bigl(V_{1}\,+\,V_{2}\bigr)\,+\,3(9\ln
N)^{2k},\le\\
\le\,\frac{3.5}{\sqrt{B}}(Bk)^{k}ML^{k-1}\,+\,0.16(Ak)^{k+1}ML^{k-1}\,=\\
=\,(Bk)^{k}ML^{k-1}\biggl(\frac{3.5}{\sqrt{B}}\,+\,0.16Ak\Bigl(\frac{e^{8}}{A}\Bigr)^{k}\biggr)\,<\,e^{9}(Bk)^{k}ML^{k-1}.
\end{multline*}
Theorem is proved.

\vspace{0.2cm}

The approximate expression for the distribution function of discrete
random quantity with the values
\[
\delta_{n}\,=\,\pi\Delta_{n}\sqrt{\frac{2}{L}},\quad N<n\le N+M,
\]
and the proof of the assertion that $\Delta_{n}\ne 0$ for `almost
all' follow from Theorem 6 by standard technic (see, for example,
Theorem 4 from \cite{Korolev_2010a}).

\vspace{0.2cm}

\begin{center}
\textbf{\S 5. On some equivalents of `almost Riemann hypothesis'}
\end{center}

The last section is devoted to some new equivalents of  `almost
Riemann hypothesis'. This hypothesis asserts that `almost all'
complex zeros of $\zeta(s)$ are on the critical line, that is
\[
\lim_{T\to +\infty}\frac{N_{0}(T)}{N(T)}\,=\,1.
\]
Moreover, the below arguments imply that Selberg interpreted Gram's
law in \cite{Selberg_1946b} in a way different from Titchmarsh's
one. Namely, the below assertions show that Selberg considered all
the complex zeros of $\zeta(s)$ (but not only the zeros on the
critical line) in dealing with the quantities $\Delta_{n}$. Thus,
the Selberg's definition of $\Delta_{n}$ is equivalent to the
Definition 4.

Suppose that $0<c_{1}<c_{2}<\ldots \le c_{n}\le c_{n+1}\le \ldots$
are the ordinates of zeros of $\zeta(s)$, lying on the critical line
and counting with theirs multiplicities. For a fixed $n\ge 1$, we
define the number $m = m(n)$ by the inequalities
\begin{equation}
t_{m-1}<c_{n}\le t_{m} \label{Lab15}
\end{equation}
and set $D_{n} = m-n$. Of course, if Riemann hypothesis is true then
$c_{n} = \gamma_{n}$ and $D_{n}=\Delta_{n}$ for any $n$.

\vspace{0.2cm}

\textsc{Theorem 7}. \emph{The validity of the relation}
\begin{equation}
\sum\limits_{n\le N}|D_{n}|\,=\,o(N^{2}) \label{Lab16}
\end{equation}
\emph{as} $N\to +\infty$, \emph{is the necessary and sufficient
condition for the truth of `almost Riemann hypothesis'}.

\vspace{0.2cm}

\textsc{Proof.} Suppose that `almost Riemann hypothesis' is true.
Then
\[
N(c_{n}+0)\,=\,(1\,+\,o(1))N_{0}(c_{n}+0)\,=\,(1\,+\,o(1))(n\,+\,O(\ln
n))\,=\,n\,+\,o(n)
\]
(by the term $O(\ln n)$, we take into account the possible
multiplicity of the zero with the ordinate $c_{n}$; by Lemma 8, this
multiplicity does not exceed in order $\ln c_{n} = O(\ln n)$). On
the other hand, by (\ref{Lab15}) we have:
\[
m-1+S(t_{m-1}+0)\,<\,N(c_{n}+0)\,\le\,m+S(t_{m}+0),
\]
and hence $D_{n} = m-n = o(n)+O(\ln m) = o(n)$. Therefore,
\[
\sum\limits_{n\le N}|D_{n}|\,=\,o(N^{2}).
\]
Suppose now that the condition (\ref{Lab16}) is satisfied. Noting
that $N(c_{n}+0)\ge n+\vep(n)$, where $\vep(n)$ is the number of
zeros of $\zeta(s)$ with the condition $0<\IIm s\le c_{n}$, $\RRe s
\ne \tfrac{1}{2}$, we get:
\[
n+\vep(n)\,\le\,N(t_{m}+0)\,=\,m+S(t_{m}+0).
\]
Hence, $0\le \vep(n)\le |D_{n}|+|S(t_{m}+0)|$. Summing this
estimation over $n\le 2N$ and applying Cauchy's inequality, we
obtain:
\[
\sum\limits_{n\le 2N}\vep(n)\,\le\,\sum\limits_{n\le
2N}|D_{n}|\,+\,\sum\limits_{\substack{n\le 2N \\ m =
m(n)}}|S(t_{m}+0)|\,\le\,\sqrt{2N}\sqrt{W}\,+\,o(N^{2}),
\]
where
\[
W\,=\,\sum\limits_{\substack{n\le 2N \\ m = m(n)}}S^{2}(t_{m}+0).
\]
Let $\mu$ be the maximum value of $m(n)$ for $n\le 2N$. Then
$t_{\mu-1}<c_{2N}\le t_{\mu}$ and hence
\[
N(c_{2N})\,\ge\,N(t_{\mu-1})\,=\,\mu\,+\,O(\ln \mu).
\]
By Lemma 8, we have for $t = c_{2N}$:
\[
N(c_{2N})\,\le\,\bigl(\tfrac{5}{2}\,-\,10^{-3}\bigr)N_{0}(c_{2N})\,=\,\bigl(\tfrac{5}{2}\,-\,10^{-3}\bigr)(2N\,+\,O(\ln
N)),
\]
and therefore $\mu<5N$. Changing the order of summation in $W$, we
obtain:
\[
W\,\le\,\sum\limits_{l\le
5N}S^{2}(t_{l}+0)\sum\limits_{\substack{n\le 2N \\ m(n)=l}}1.
\]
For a fixed $l$, the number of $n$ that satisfy the conditions $n\le
2N$, $m(n) = l$, does not exceed the number of all ordinates of
zeros of $\zeta(s)$ in the interval $(t_{l-1},t_{l}]$, that is
\[
N(t_{l}+0)\,-\,N(t_{l-1}+0)\,=\,1\,+\,S(t_{l}+0)\,-\,S(t_{l
-1}+0)\,=\,1\,+\,r(l).
\]
Thus we have:
\[
W\,\le\,\sum\limits_{l\le 5N}S^{2}(t_{l}+0)(1\,+\,r(l)).
\]
Using both the first formula of Lemma 5 and the estimation of
Theorem 5, we find that
\begin{align*}
& W\,\le\,\frac{5N}{2\pi^{2}}\,\ln\ln N\,+\,O\bigl(N\sqrt{\ln\ln
N}\bigr)\,<\,\tfrac{1}{3}N\ln\ln N,\\
& \sum\limits_{n\le 2N}\vep(n)\,<\,N\sqrt{\ln\ln
N}\,+\,o(N^{2})\,=\,o(N^{2}).
\end{align*}
By obvious inequality $\vep(N+1)+\vep(N+2)+\ldots + \vep(2N)\ge
N\vep(N)$ we get:
\[
N\vep(N)\,=\,o(N^{2}),\quad \vep(N)\,=\,o(N).
\]
Suppose now $t$ is a sufficiently large. Then, defining $N$ from the
inequalities $c_{N-1}<t\le c_{N}$ and using the above relations, we
obtain:
\[
N(t)\,\ge\, N-1,\quad
N(t)\,-\,N_{0}(t)\,\le\,\vep(N)\,=\,o(N)\,=\,o(N(t)).
\]
The theorem is proved.

\vspace{0.2cm}

\textsc{Coollary 1.} \emph{The validity of the relation}
\[
\sum\limits_{n\le N}|D_{n}|^{k}\,=\,o\bigl(N^{k+1}\bigr),\quad N\to
+\infty
\]
\emph{for at least one fixed value of} $k\ge 1$ \emph{is the
necessary and sufficient condition for the truth of the `almost
Riemann hypothesis'}

\vspace{0.2cm}

The proof is similar to the previous one. The difference is that we
should use the inequality
\[
N\vep(N)\,\le\,\sum\limits_{n\le 2N}|D_{n}|\,+\,N\sqrt{\ln\ln
N}\,\le\,(2N)^{1-1/k}\biggl(\,\sum\limits_{n\le
2N}|D_{n}|^{k}\biggr)^{\!\!1/k}\,+\,N\sqrt{\ln\ln N}
\]
for the proof of sufficiency.

This assertion shows, in particular, that if Selberg's formulas
(\ref{Lab08}),(\ref{Lab09}) hold true after the replacement of the
quantities $\Delta_{n}$ by $D_{n}$, then the `almost Riemann
hypothesis' is also true.

\vspace{0.2cm}

\textsc{Corollary 2}. \emph{The assertion} `$D_{n} = o(n)$ \emph{as}
$n\to +\infty$' \emph{is the necessary and sufficient condition for
the `almost Riemann hypothesis'}.

\vspace{0.2cm}

The below theorem shows that the upper bound for $D_{n}$ causes the
main difficulty.

\vspace{0.2cm}

\textsc{Theorem 8}. \emph{Suppose that} $N_{0}(t)>\vk N(t)$
\emph{for any} $t>t_{0}>1$ \emph{and for some constant} $\vk$,
$0<\vk<1$. \emph{Then the following inequalities hold for all
sufficiently large} $n$:
\[
-9\ln{n}\,\le\,D_{n}\,\le\,\Bigl(\frac{1}{\vk}\,-\,1\Bigr)n\,+\,9\Bigl(\frac{1}{\vk}\,+\,1\Bigr)\ln{n}.
\]

\vspace{0.2cm}

\textsc{Proof.} By (\ref{Lab15}), we get:
\begin{equation}
N(c_{n}+0)\,\ge\,N(t_{m-1}\,+\,0)\,=\,m\,-\,1\,+\,S(t_{m-1}+0).\label{Lab17}
\end{equation}
By the assumption of theorem, we get for $t = c_{n}+0$:
\begin{equation}
N(c_{n}+0)\,<\,\frac{1}{\vk}\,N_{0}(c_{n}+0)\,\le\,\frac{1}{\vk}\,(n\,+\,\kappa_{n}),\label{Lab18}
\end{equation}
where $\kappa_{n}$ denotes the multiplicity of the ordinate $c_{n}$.
Comparing (\ref{Lab17}) and (\ref{Lab18}) and using the inequality
$|S(t)|\le 8.9\ln t$, we obtain:
\[
m\,\le\,\frac{1}{\vk}\,n\,+\,9\Bigl(\frac{1}{\vk}\,+\,1\Bigr)\ln{n},
\quad
D_{n}\,=\,m-n\,<\,\Bigl(\frac{1}{\vk}\,-\,1\Bigr)n\,+\,9\Bigl(\frac{1}{\vk}\,+\,1\Bigr).
\]
On the other hand,
\[
n\,\le\,N_{0}(c_{n}+0)\,\le\,N(c_{n}+0)\,\le\,N(t_{m}+0)\,=\,m\,+\,S(t_{m}+0),
\]
and therefore $D_{n}\ge -S(t_{m}+0)\ge -9\ln{n}$. Theorem is proved.

\vspace{0.2cm}

In \cite{Selberg_1946b}, Selberg referred to the formulas
\begin{equation}
\liminf_{n\to +\infty}\Delta_{n}\,=\,-\infty,\quad \limsup_{n\to
+\infty}\Delta_{n}\,=\,+\infty,\label{Lab19}
\end{equation}
as Tithmarsh's result from \cite{Titchmarsh_1935}. Indeed, the
relations (\ref{Lab19}) hold true, and the modern omega\,-theorems
for the function $S(t)$ imply a much deeper result, namely
\[
\Delta_{n}\,=\,\Omega_{\pm}\biggl(\sqrt[3\;]{\frac{\ln{n}}{\ln\ln{n}}}\biggr),
\]
as $n$ growths (see \cite{Korolev_2010b}). On the other hand, in
\cite{Titchmarsh_1935}, Titchmarsh considered the fractions
\[
\tau_{n}\,=\,\frac{c_{n}\,-\,t_{n}}{t_{n+1}\,-\,t_{n}}
\]
instead of the quantities $\Delta_{n}$ (one can easily see that the
difference between $\tau_{n}$ and $D_{n}$ is $O(1)$), and
established the unboundedness of $\tau_{n}$. As far as can be seen,
the methods of \cite{Titchmarsh_1935} allows one to show only that
$\tau_{n}\ne O(1)$ and $D_{n}\ne O(1)$, as $n\to +\infty$. Slight
modification of these methods and the omega\,-theorems for $S(t)$
lead to the following assertion
\[
D_{n}\,=\,\Omega_{-}\biggl(\sqrt[3\;]{\frac{\ln{n}}{\ln\ln{n}}}\biggr),\quad
n\to +\infty.
\]
So, the problem of unboundedness of $D_{n}$ from above remains still
open.

\renewcommand{\refname}{

\end{document}